\documentclass[11pt,a4paper]{article}
\usepackage[margin=28mm,headheight=14pt]{geometry}
\usepackage{amsmath,amssymb,amsthm,mathtools}
\usepackage{lmodern}
\usepackage{microtype}
\usepackage{tikz}
\usepackage{enumitem}
\usepackage{float}
\usepackage{placeins}
\usepackage{hyperref}
\hypersetup{hidelinks,pdftitle={Zeros of Pick and Hopf Differentials and Fences of Convex Domains},pdfauthor={}}
\usepackage{fancyhdr}
\allowdisplaybreaks[2]
\newtheorem{mainthm}{Theorem}

\newtheorem{theorem}{Theorem}[section]
\newtheorem{lemma}[theorem]{Lemma}
\newtheorem{proposition}[theorem]{Proposition}
\newtheorem{corollary}[theorem]{Corollary}
\theoremstyle{definition}

\theoremstyle{remark}
\newtheorem{remark}[theorem]{Remark}
\newtheorem*{remark*}{Remark}
\numberwithin{equation}{section}

\DeclareMathOperator{\tr}{tr}
\DeclareMathOperator{\Int}{int}

\newcommand{\R}{\mathbb R}
\newcommand{\C}{\mathbb C}
\newcommand{\HH}{\mathbb H^2}
\newcommand{\RP}{\mathbb{RP}^2}
\newcommand{\Id}{\mathrm{Id}}

\newcommand{\II}{\mathrm{II}}
\newcommand{\Rea}{\operatorname{Re}}

\definecolor{fencefill}{RGB}{218,232,239}
\definecolor{fenceink}{RGB}{42,77,98}
\definecolor{fencefilldeep}{RGB}{161,193,210}

\title{\bfseries Zeros of Pick and Hopf Differentials and Fences of Convex Domains}
\author{Qiongling Li\footnote{qiongling.li@nankai.edu.cn, Chern Institute of Mathematics and LPMC, Nankai University}, 
Erxiao Wang\footnote{wang.eric@zjnu.edu.cn, Zhejiang Normal University}} 
\date{}
\date{}
\makeatletter
\renewcommand{\maketitle}{%
  \thispagestyle{plain}%
  \begin{center}
    {\LARGE\bfseries\@title\par}%
    \ifx\@author\@empty\else\vspace{.7em}{\normalsize\@author\par}\fi
  \end{center}\vspace{.5em}}
\makeatother
\begin{document}
\maketitle
\begin{abstract}
We prove that the vanishing order of a Pick or Hopf differential
controls the corresponding convex domain through sharp disk
comparisons. For a complete hyperbolic affine sphere, Pick
vanishing order $m\ge0$ gives normalized projective coordinates
in which the domain lies between disks of radii
$\sqrt{\cos(\pi/(m+3))}$ and its reciprocal. For a harmonic
diffeomorphism to a convex hyperbolic domain with complete
$\partial$-energy metric, Hopf vanishing order $m$ gives a disk
in the image of radius $\log\cot(\pi/(2(m+2)))$, centered at
the image point. Equality characterizes the corresponding regular
polygons. For polygonal targets, these comparisons prove the
Fence conjecture of Dumas--Wolf and give explicit Hilbert and
hyperbolic distance bounds for higher fences. We also construct simple-zero examples showing that these
sufficient multiplicity conditions are not necessary.


\vskip 5pt
{\em 2020 MR Subject Classification}: Primary 53A15; Secondary 35J96, 52A10.
\end{abstract}

\section{Introduction}
Holomorphic cubic and quadratic differentials provide analytic
descriptions of convex domains in the projective and hyperbolic
planes. We study how the vanishing order of a holomorphic differential
at a single point constrains the geometry of its associated
convex domain. 

On the cubic side, the Cheng--Yau theorem associates to a
properly convex domain $\Omega\subset\RP$ a complete hyperbolic
affine sphere with Blaschke metric $h$ and holomorphic Pick
differential $C$ \cite{CY}. Radial projection identifies the sphere
with $\Omega$. Wang and Simon--Wang described its reconstruction
from a cubic differential and a conformal metric satisfying Wang's
equation \cite{Wang91,SimonWang}. Dumas and Wolf established a
homeomorphism between the moduli spaces of polynomial cubic
differentials of degree $n-3$ and oriented convex projective
$n$-gons \cite{DW}. Related correspondences treat
bounded cubic differentials on the disk \cite{BH} and equivariant
data on closed and punctured surfaces
\cite{LabourieCubic,LoftinAffine,BHFinite,Nie}.

The quadratic counterpart is an orientation-preserving harmonic
diffeomorphism $f:X\to\Omega$ onto a convex domain
$\Omega\subset\HH$, with holomorphic Hopf differential $Q_f$. Such a map can be realized as the Gauss map of
a strictly convex spacelike constant mean curvature (CMC) $1$ surface
in Minkowski three-space. The induced metric is the $\partial$-energy metric $I_f$ of $f$, and
Hopf zeros are precisely the umbilic points. Wan and Wan--Au
developed the existence theory for prescribed Hopf differentials
through this correspondence \cite{Wan,WanAu}, while Choi and
Treibergs studied the geometry of the Gauss images
\cite{ChoiTreibergs}. Bonsante, Seppi and Smillie developed the
asymptotic theory of entire constant-curvature surfaces and CMC
hypersurfaces in Minkowski space \cite{BSS19,BSS}.
Under completeness of the $\partial$-energy metric, Han, Tam,
Treibergs and Wan characterized polynomial Hopf
differentials on $\C$ of degree $n-2\ge1$ by ideal $n$-gonal images
\cite{HTTW}. Au, Tam and Wan studied nonpolynomial data and their
image boundaries \cite{ATW}, and Au and Wan described image
polygons for symmetric polynomial data \cite{AuWanSymmetry}.
For bounded quadratic differentials on the disk, the correspondence
yields harmonic quasiconformal self-diffeomorphisms \cite{Wan}.
Equivariant versions connect quadratic differentials with
Teichm\"uller spaces of closed, punctured, and crowned surfaces
\cite{Wolf,Lohkamp,Gupta}.

Both constructions are governed by scalar vortex equations,
the rank-two and rank-three cases of cyclic Toda systems.
Existence and uniqueness results for complete solutions were studied in \cite{LiVortex, LiMochizuki}; in these two cases, completeness refers to
the $\partial$-energy and Blaschke metrics, respectively.
Numerical studies of both polynomial correspondences appear
in \cite{DN}.

\medskip

Our main results are a pair of sharp disk comparisons
depending only on the vanishing order, with regular polygons
as the equality models.  We first
state these polygonal consequences, including the Fence conjecture
of Dumas--Wolf, and then the disk comparisons for general domains.

\subsection{Higher fences}
The Fence conjecture of Dumas and Wolf asserts that all Pick zeros
lie in the region cut out by the short diagonals joining the two
neighbors of each vertex \cite[Conjecture~9.2]{DW}.
We call this region the \textbf{ordinary fence} and also consider higher
fences defined by diagonals that skip more vertices. 

\begin{figure}[H]
\centering
\begin{tikzpicture}[scale=2.35,line cap=round,line join=round]
\path[use as bounding box] (-1.22,-1.66) rectangle (4.32,1.07);
\tikzset{
  fencepoly/.style={draw=black!75,line width=.65pt},
  fencetwo/.style={draw=fenceink!75,line width=.75pt,dash pattern=on 2.2pt off 1.25pt},
  fencethree/.style={draw=fenceink,line width=.9pt},
  fencecutlight/.style={draw=fenceink!42,line width=.38pt,dash pattern=on 1.2pt off 1.1pt},
  fencecutdeep/.style={draw=fenceink!55,line width=.38pt,densely dotted}}
\begin{scope}[shift={(0.0,0)}]
\fill[black!3] (1.0000000000,0.0464993911) -- (0.9441017801,0.3502787946) -- (0.2289818182,0.8538147727) -- (-0.8658468566,0.5264820659) -- (-1.0000000000,-0.4690059819) -- (-0.1800535995,-0.8538147727) -- (0.6591975999,-0.5748602758) -- cycle;
\fill[fencefill] (0.8869908801,0.1648286379) -- (0.6858693055,0.3754184068) -- (-0.1407160028,0.4558886656) -- (-0.5903285088,-0.0280540141) -- (-0.3506611254,-0.5104328186) -- (0.2221244421,-0.5469756703) -- (0.8046526330,-0.1025396030) -- cycle;
\fill[fencefilldeep] (0.4440667021,0.1395536508) -- (0.0145328304,-0.0414608236) -- (0.0006789862,-0.0992974625) -- (0.1147499737,-0.1816761814) -- (0.5101116533,-0.0797706464) -- cycle;
\draw[fencecutlight] (1.0000000000,0.0464993911) -- (0.2289818182,0.8538147727);
\draw[fencecutlight] (0.9441017801,0.3502787946) -- (-0.8658468566,0.5264820659);
\draw[fencecutlight] (0.2289818182,0.8538147727) -- (-1.0000000000,-0.4690059819);
\draw[fencecutlight] (-0.8658468566,0.5264820659) -- (-0.1800535995,-0.8538147727);
\draw[fencecutlight] (-1.0000000000,-0.4690059819) -- (0.6591975999,-0.5748602758);
\draw[fencecutlight] (-0.1800535995,-0.8538147727) -- (1.0000000000,0.0464993911);
\draw[fencecutlight] (0.6591975999,-0.5748602758) -- (0.9441017801,0.3502787946);
\draw[fencecutdeep] (1.0000000000,0.0464993911) -- (-0.8658468566,0.5264820659);
\draw[fencecutdeep] (0.9441017801,0.3502787946) -- (-1.0000000000,-0.4690059819);
\draw[fencecutdeep] (0.2289818182,0.8538147727) -- (-0.1800535995,-0.8538147727);
\draw[fencecutdeep] (-0.8658468566,0.5264820659) -- (0.6591975999,-0.5748602758);
\draw[fencecutdeep] (-1.0000000000,-0.4690059819) -- (1.0000000000,0.0464993911);
\draw[fencecutdeep] (-0.1800535995,-0.8538147727) -- (0.9441017801,0.3502787946);
\draw[fencecutdeep] (0.6591975999,-0.5748602758) -- (0.2289818182,0.8538147727);
\draw[fencetwo] (0.8869908801,0.1648286379) -- (0.6858693055,0.3754184068) -- (-0.1407160028,0.4558886656) -- (-0.5903285088,-0.0280540141) -- (-0.3506611254,-0.5104328186) -- (0.2221244421,-0.5469756703) -- (0.8046526330,-0.1025396030) -- cycle;
\draw[fencethree] (0.4440667021,0.1395536508) -- (0.0145328304,-0.0414608236) -- (0.0006789862,-0.0992974625) -- (0.1147499737,-0.1816761814) -- (0.5101116533,-0.0797706464) -- cycle;
\draw[fencepoly] (1.0000000000,0.0464993911) -- (0.9441017801,0.3502787946) -- (0.2289818182,0.8538147727) -- (-0.8658468566,0.5264820659) -- (-1.0000000000,-0.4690059819) -- (-0.1800535995,-0.8538147727) -- (0.6591975999,-0.5748602758) -- cycle;
\fill[black!85] (1.0000000000,0.0464993911) circle[radius=.010];
\fill[black!85] (0.9441017801,0.3502787946) circle[radius=.010];
\fill[black!85] (0.2289818182,0.8538147727) circle[radius=.010];
\fill[black!85] (-0.8658468566,0.5264820659) circle[radius=.010];
\fill[black!85] (-1.0000000000,-0.4690059819) circle[radius=.010];
\fill[black!85] (-0.1800535995,-0.8538147727) circle[radius=.010];
\fill[black!85] (0.6591975999,-0.5748602758) circle[radius=.010];
\fill[white] (0.1151975983,-0.0792158299) circle[radius=.024];
\fill[black] (0.1151975983,-0.0792158299) circle[radius=.014];
\node[font=\small,anchor=west,inner sep=0pt] at (0.1801975983,-0.0392158299) {$p$};
\node[font=\small] at (0,-1.18) {(a) Projective polygon: $m=3$};
\end{scope}
\begin{scope}[shift={(3.1,0)}]
\fill[black!3] (0.9899941743,-0.1411082380) arc[start angle=-98.1119802673,delta angle=-163.5420313245,radius=0.1446186158] arc[start angle=-81.6540115917,delta angle=-132.2270676880,radius=0.4428564178] arc[start angle=-33.8810822730,delta angle=-102.5876089703,radius=0.8013286067] arc[start angle=43.5313087567,delta angle=-93.0783722590,radius=0.9476655833] arc[start angle=130.4529364977,delta angle=-122.1349195419,radius=0.5527708948] arc[start angle=-171.6819830442,delta angle=-138.5338866356,radius=0.3785279669] arc[start angle=-130.2158733729,delta angle=-147.8961068943,radius=0.2877268062] -- cycle;
\fill[fencefill] (0.7011442095,-0.0212412440) arc[start angle=-126.9629767020,delta angle=-19.1965386612,radius=0.6276746326] arc[start angle=-94.5278690307,delta angle=-13.8338553834,radius=1.9285908358] arc[start angle=-39.8674895657,delta angle=-2.9404926770,radius=7.2690469832] arc[start angle=26.1381874305,delta angle=-6.1339575107,radius=3.1511311964] arc[start angle=90.7539022998,delta angle=-15.6779748455,radius=1.1777251887] arc[start angle=140.6265663953,delta angle=-29.8928441477,radius=0.7476873820] arc[start angle=164.6867159135,delta angle=-23.0441720862,radius=0.4511166814] -- cycle;
\fill[fencefilldeep] (0.3870109877,0.0009838815) arc[start angle=108.5718744744,delta angle=3.7591916393,radius=3.4751701672] arc[start angle=168.0664781339,delta angle=0.8696980378,radius=2.5916213021] arc[start angle=-132.7136455276,delta angle=0.2288836408,radius=18.3081883459] arc[start angle=101.9196878011,delta angle=-5.2857111343,radius=2.2165636355] arc[start angle=-162.2176029819,delta angle=-7.1340858998,radius=1.1171590402] -- cycle;
\draw[fencecutlight] (0.9899941743,-0.1411082380) arc[start angle=-98.1119802673,delta angle=-115.7690995801,radius=0.6276746326];
\draw[fencecutlight] (0.9894096030,0.1451503964) arc[start angle=-81.6540115917,delta angle=-54.8146796516,radius=1.9285908358];
\draw[fencecutlight] (0.5574710042,0.8301963922) arc[start angle=-33.8810813250,delta angle=-15.6659821773,radius=7.2690469832];
\draw[fencecutlight] (-0.6887508471,0.7249981176) arc[start angle=43.5313087567,delta angle=-35.2132918009,radius=3.1511311964];
\draw[fencecutlight] (-0.7609391741,-0.6488232219) arc[start angle=130.4529364977,delta angle=-80.6688072598,radius=1.1777251887];
\draw[fencecutlight] (0.1446673554,-0.9894803466) arc[start angle=-171.6819830442,delta angle=-106.4299972231,radius=0.7476873820];
\draw[fencecutlight] (0.7636171956,-0.6456692322) arc[start angle=-130.2158726130,delta angle=-131.4381389787,radius=0.4511166814];
\draw[fencecutdeep] (0.9899941743,-0.1411082380) arc[start angle=-98.1119802673,delta angle=-38.3567109761,radius=2.8751055290];
\draw[fencecutdeep] (0.9894096030,0.1451503964) arc[start angle=98.3459884083,delta angle=32.1069480894,radius=3.4751701672];
\draw[fencecutdeep] (0.5574710042,0.8301963922) arc[start angle=146.1189191219,delta angle=42.1990978339,radius=2.5916213021];
\draw[fencecutdeep] (-0.6887508471,0.7249981176) arc[start angle=-136.4686912433,delta angle=6.2528199356,radius=18.3081883459];
\draw[fencecutdeep] (-0.7609391741,-0.6488232219) arc[start angle=130.4529364977,delta angle=-48.5649167650,radius=2.2165636355];
\draw[fencecutdeep] (0.1446673554,-0.9894803466) arc[start angle=-171.6819830442,delta angle=-89.9720285475,radius=1.0004883131];
\draw[fencecutdeep] (0.7636171956,-0.6456692322) arc[start angle=-130.2158718151,delta angle=-83.6652086282,radius=1.1171590402];
\draw[fencetwo] (0.7011442095,-0.0212412440) arc[start angle=-126.9629767020,delta angle=-19.1965386612,radius=0.6276746326] arc[start angle=-94.5278690307,delta angle=-13.8338553834,radius=1.9285908358] arc[start angle=-39.8674895657,delta angle=-2.9404926770,radius=7.2690469832] arc[start angle=26.1381874305,delta angle=-6.1339575107,radius=3.1511311964] arc[start angle=90.7539022998,delta angle=-15.6779748455,radius=1.1777251887] arc[start angle=140.6265663953,delta angle=-29.8928441477,radius=0.7476873820] arc[start angle=164.6867159135,delta angle=-23.0441720862,radius=0.4511166814] -- cycle;
\draw[fencethree] (0.3870109877,0.0009838815) arc[start angle=108.5718744744,delta angle=3.7591916393,radius=3.4751701672] arc[start angle=168.0664781339,delta angle=0.8696980378,radius=2.5916213021] arc[start angle=-132.7136455276,delta angle=0.2288836408,radius=18.3081883459] arc[start angle=101.9196878011,delta angle=-5.2857111343,radius=2.2165636355] arc[start angle=-162.2176029819,delta angle=-7.1340858998,radius=1.1171590402] -- cycle;
\draw[fencepoly] (0.9899941743,-0.1411082380) arc[start angle=-98.1119802673,delta angle=-163.5420313245,radius=0.1446186158] arc[start angle=-81.6540115917,delta angle=-132.2270676880,radius=0.4428564178] arc[start angle=-33.8810822730,delta angle=-102.5876089703,radius=0.8013286067] arc[start angle=43.5313087567,delta angle=-93.0783722590,radius=0.9476655833] arc[start angle=130.4529364977,delta angle=-122.1349195419,radius=0.5527708948] arc[start angle=-171.6819830442,delta angle=-138.5338866356,radius=0.3785279669] arc[start angle=-130.2158733729,delta angle=-147.8961068943,radius=0.2877268062] -- cycle;
\draw[black!30,line width=.35pt] (0,0) circle[radius=1];
\draw[black!75,fill=white,line width=.35pt] (0.9899941743,-0.1411082380) circle[radius=.010];
\draw[black!75,fill=white,line width=.35pt] (0.9894096030,0.1451503964) circle[radius=.010];
\draw[black!75,fill=white,line width=.35pt] (0.5574710125,0.8301964046) circle[radius=.010];
\draw[black!75,fill=white,line width=.35pt] (-0.6887508471,0.7249981176) circle[radius=.010];
\draw[black!75,fill=white,line width=.35pt] (-0.7609391741,-0.6488232219) circle[radius=.010];
\draw[black!75,fill=white,line width=.35pt] (0.1446673554,-0.9894803466) circle[radius=.010];
\draw[black!75,fill=white,line width=.35pt] (0.7636172037,-0.6456692390) circle[radius=.010];
\fill[white] (0.2200000000,-0.1000000000) circle[radius=.024];
\fill[black] (0.2200000000,-0.1000000000) circle[radius=.014];
\node[font=\small,anchor=west,inner sep=0pt] at (0.2850000000,-0.0600000000) {$p$};
\node[font=\small] at (0,-1.18) {(b) Ideal polygon: $m=4$};
\end{scope}
\begin{scope}[shift={(1.55,-1.51)}]
\draw[fencepoly,fill=black!3] (-1.18,-.055) rectangle (-.94,.055);
\node[anchor=west,font=\small,inner sep=3pt] at (-.94,0) {$P$};
\draw[fencetwo,fill=fencefill] (-.38,-.055) rectangle (-.14,.055);
\node[anchor=west,font=\small,inner sep=3pt] at (-.14,0) {$F_2(P)$};
\draw[fencethree,fill=fencefilldeep] (.62,-.055) rectangle (.86,.055);
\node[anchor=west,font=\small,inner sep=3pt] at (.86,0) {$F_3(P)$};
\end{scope}
\end{tikzpicture}
\caption{\small Higher fences of an irregular projective heptagon in an affine chart
and an irregular ideal heptagon in the Poincar\'e disk, with all defining
diagonals drawn. The markers schematically illustrate Theorem~\ref{thm:higher-fence-intro}
for a Pick zero of order $3$ in (a) and the image of a Hopf zero of order $4$
in (b). In both cases, $p\in F_3(P)\subset F_2(P)$.}
\label{fig:higher-fences}
\end{figure}
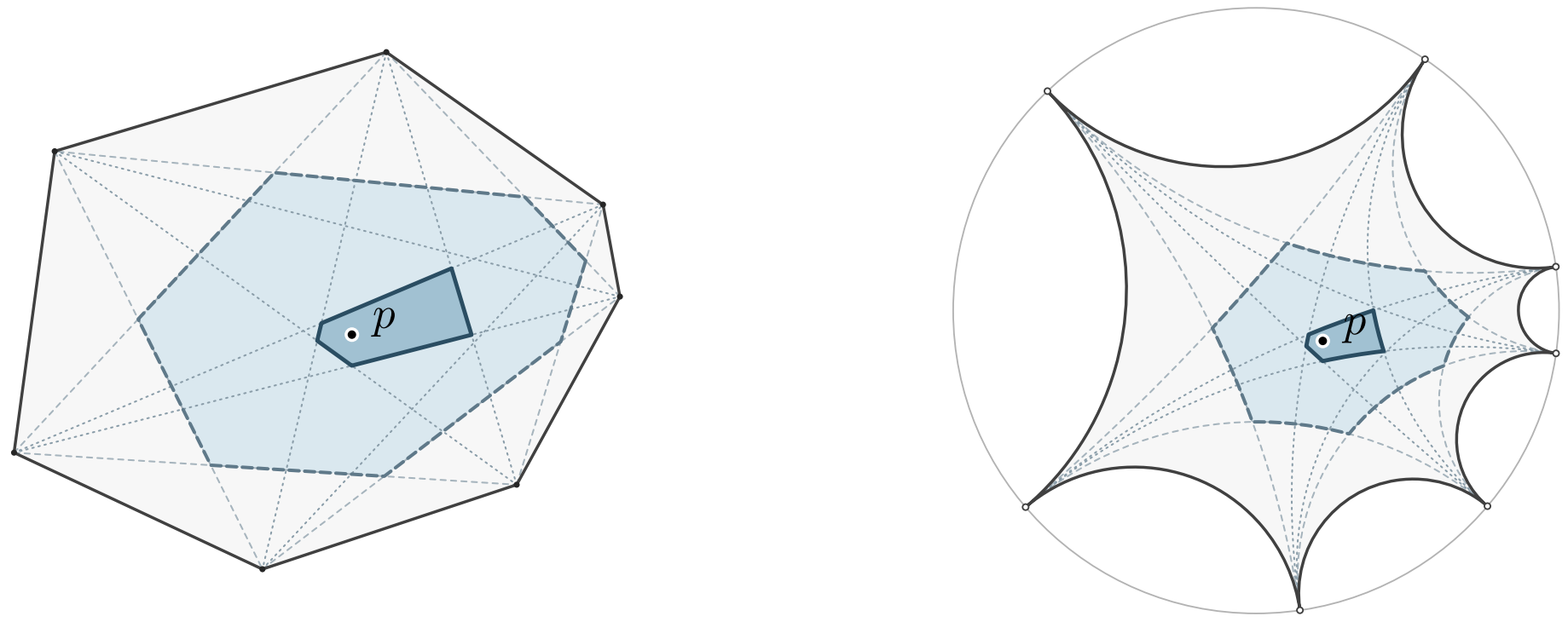

Let $P$ be a convex projective polygon or an ideal hyperbolic polygon,
with distinct vertices $v_1,\ldots,v_n$ in cyclic order and indices
taken modulo $n$. Vertices in the projective case are genuine corners.
For $2\le s\le n-1$, define
\begin{equation}\label{eq:fence-def}
 F_s(P)=\overline P\cap\bigcap_{i=1}^n H_i^{(s)},
\end{equation}
where $H_i^{(s)}$ is the closed half-plane bounded by the line through
$v_i,v_{i+s}$ and excluding the intermediate vertices
$v_{i+1},\ldots,v_{i+s-1}$.
In the hyperbolic case, the bounding line is the geodesic with ideal
endpoints $v_i,v_{i+s}$; take the closure in the closed-disk
compactification and then intersect with $\HH$. We use Klein coordinates in the proofs. Figures 1(b) and 2(b) are drawn in the Poincaré disk.
Thus $F_2(P)$ is the ordinary fence, while larger $s$ gives higher
fences by skipping $s-1$ vertices. These regions are empty when
$2s>n$; see Section~\ref{subsec:polygon-fences}.
Figure~\ref{fig:higher-fences} illustrates the nested regions
$F_3(P)\subset F_2(P)\subset P$ in the projective and
ideal-hyperbolic settings. 

\FloatBarrier

The vanishing order gives a uniform criterion for membership in
higher fences, together with an explicit lower bound for the
distance to their boundary.

\begin{mainthm}\label{thm:higher-fence-intro}
Let $P$ have $n$ vertices and let $m\ge0$. Consider either of the
following settings:
\begin{enumerate}[label=\textup{(\roman*)},leftmargin=2em,beginpenalty=10000]
\item $P$ is a convex projective polygon and its Pick differential has
vanishing order $m$ at $p\in\Int P$. Set $k=m+3$ and let $d=d_P$ be
the Hilbert metric.
\item $f:X\to P$ is an orientation-preserving harmonic diffeomorphism
onto an ideal polygon with complete $\partial$-energy metric, and
the Hopf differential has vanishing order $m$ at $z_0\in X$.
Set $p=f(z_0)$, $k=m+2$, and $d=d_{\HH}$.
\end{enumerate}
For every $2\le s\le n-1$ with $k\ge2s$, we have $p\in F_s(P)$ and
\[
 d\bigl(p,\partial F_s(P)\bigr)
 \ge\log\cot\frac{s\pi}{2k}.
\]
\end{mainthm}

 The polynomial degrees give $k\le n$ in both settings.
Corollary~\ref{cor:higher-fences} proves the theorem and gives its
critical and equality cases. Corollary~\ref{cor:general-fences}
extends the bounds to general convex domains using finite polygonal
boundary chains. The multiplicity conditions are sufficient rather
than necessary: Section~\ref{sec:examples} constructs polygons
whose zero images all lie in a prescribed higher fence even though
every zero is simple.

\FloatBarrier
Taking $s=2$ recovers the ordinary fence, including the critical
case described in Corollary~\ref{cor:higher-fences}.

\begin{corollary}\label{cor:fence-intro}
In the setting of Theorem~\ref{thm:higher-fence-intro}, suppose
$n\ge4$ and assume $m\ge1$ in the Pick case or $m\ge2$ in the Hopf case.
Then $p\in\Int F_2(P)$ unless $n=4$, in which case
$F_2(P)=\{p\}$.
\end{corollary}

The Pick case proves the Fence conjecture of Dumas and Wolf.
The multiplicity restriction in the Hopf case is essential:
the polynomial quadratic differential $(z^2-1)\,dz^2$ produces a
harmonic diffeomorphism onto an ideal quadrilateral
\cite[Theorem~1.1 and \S2]{HTTW}. Its two distinct zero images
cannot both equal the single point of that quadrilateral's fence.

\subsection{Estimates for general convex domains}
The disk comparisons behind the higher-fence theorem apply
without requiring a polygonal boundary or a prescribed conformal type of
the domain.

For the projective comparison, a chart is normalized at $p$ if
$p=0$ and its canonical affine-sphere potential $U$, recalled in
\eqref{eq:affine-pde}, satisfies
\[
 U(0)=-1,\qquad DU(0)=0,\qquad D^2U(0)=\Id.
\]
Such a chart exists by Lemma~\ref{lem:affine-contact} and is unique
up to an orthogonal linear change of coordinates. The centered
radii and regularity of the normalized polygon are therefore
independent of this choice.

\begin{samepage}
\begin{mainthm}\label{thm:pick}
Let $\Omega\subset\RP$ be properly convex, with complete hyperbolic
affine sphere $S_\Omega$, Blaschke metric $h$, and Pick differential
$C$. Suppose $C$ has vanishing order $m\ge0$ over $p\in\Omega$,
and put $k=m+3$.
There are bounded projective affine coordinates normalized at $p$
with $h_p=|dx|^2$, such that
\begin{equation}\label{eq:pick-main}
 B_{r_k}(0)\Subset\Omega\Subset B_{r_k^{-1}}(0),
 \qquad r_k=\sqrt{\cos(\pi/k)},
\end{equation}
unless the normalized domain is a regular $k$-gon. In that case,
its inradius and circumradius about $0$ are $r_k$ and $r_k^{-1}$,
respectively, so both radius bounds are attained.
\end{mainthm}
\end{samepage}

In particular, a high-order Pick zero forces the normalized domain
to be nearly circular. Corollary~\ref{cor:Hilbert-norm} also gives
vanishing-order-dependent comparisons between the Hilbert and
Blaschke norms, with multiplicative error $1+O(m^{-2})$ as
$m\to\infty$. The corresponding Hopf estimate is intrinsic to
the hyperbolic target.

\begin{mainthm}\label{thm:hopf}
Let $f:X\to\Omega\subset\HH$ be an orientation-preserving harmonic
diffeomorphism onto a convex open domain. Assume that $I_f$ is
complete. If $Q_f$ has vanishing order $m\ge0$ at $z_0$, then
\begin{equation}\label{eq:hopf-main}
 B_{\HH}\bigl(f(z_0),r_m\bigr)\Subset\Omega,
 \qquad r_m=\log\cot\frac{\pi}{2(m+2)},
\end{equation}
unless $m\ge1$ and $\Omega$ is a regular ideal $(m+2)$-gon
centered at $f(z_0)$. In that case,
$d_{\HH}(f(z_0),\partial_{\HH}\Omega)=r_m$.
\end{mainthm}

Under the completeness assumption, the image is the interior of
the convex hull of its ideal boundary; see \eqref{eq:support-domain}.
It need not have finitely many sides.

Thus $d_{\HH}(f(z_0),\partial_{\HH}\Omega)\ge r_m$, with
equality precisely in the regular ideal polygon case just described.
The boundary is taken inside $\HH$, and the distance is $+\infty$
when $\Omega=\HH$. For $m=0$, the distance is positive at every
interior point, but its best uniform lower bound is $r_0=0$.
Already at a simple Hopf zero, the theorem gives the positive bound
\[
 d_{\HH}\bigl(f(z_0),\partial_{\HH}\Omega\bigr)
 \ge \tfrac12\log3,
\]
with equality exactly for a regular ideal triangle centered at
$f(z_0)$.

\begin{remark}
For $m\ge1$, equality in Theorem~\ref{thm:hopf} also determines
the harmonic map. Indeed, the proof in Section~\ref{sec:hopf}
shows that, after a suitable normalization, its CMC realization
coincides with the regular $(m+2)$-gonal model.
Consequently, $f$ is equivalent to the Gauss map of this model,
up to an orientation-preserving isometry of $\HH$ and a conformal
reparametrization of the source.
\end{remark}

\subsection{Method of proof}
The proofs compare the affine-sphere and CMC realizations with
regular polygonal models. For the nontrivial disk comparisons,
normalization at a point of vanishing order $m$ makes the defining
potential or graph function
agree with a quadric up to an error $O(|x|^{m+r})$, where $r=3$
for Pick and $r=2$ for Hopf differentials. The corresponding
regular $(m+r)$-gonal model has the same estimate. Thus, unless
the surfaces coincide, their difference vanishes to order at
least $m+r$. Hartman--Wintner local nodal theory
\cite{HartmanWintner,Alessandrini}, a planar separation argument
and the maximum principle then give at least $m+r$ negative
components. The global step is to control these components
using the geometry of the model.

For affine spheres, the components meet distinct sides of the
model polygon. Rotating the model yields the inner disk bound.
Conormal polarity determines the radii of the normalized regular
models and converts the inner bound into the outer one.

For CMC graphs, the comparison takes place over the entire plane,
so there is no finite boundary on which to count components
directly. We pass to Legendre transforms and use the
support-function theory of Bonsante, Seppi and Smillie
\cite{BSS19,BSS}. On the dual polygon the model support function is finite, whereas the support function of the surface being compared may take the value $+\infty$. A convex patching argument associates every
negative comparison component with a vertex of the model's Gauss
image. If the other Gauss image misses a model edge, the upper bound on the number of
components therefore drops by one. Comparing this bound with
the local nodal count proves the disk inclusion and its rigidity.

A common convex-geometric estimate then converts the disk
comparisons into the higher-fence bounds.
\subsection{Related work}
Spacelike maximal surfaces in $\mathbb H^{2,n}$ exhibit analogous
links between holomorphic data and asymptotic boundaries.
In $\mathrm{AdS}^3$, these include minimal Lagrangian extensions
of quasisymmetric maps \cite{BSMaximal} and polynomial
parametrizations of lightlike polygons \cite{TamburelliPolygons}.
Higher-codimension developments include equivariant maximal
surfaces for maximal representations \cite{CTT}, the asymptotic
Plateau problem \cite{LTW}, and the correspondence between
quasisymmetric curves and quasiperiodic maximal surfaces
\cite{LTQuasicircles}. The relation between polynomial data and polygonal maximal
surfaces in higher codimension has been further developed in
\cite{TW,MorianiPolygons,EvansMoriani}. It is natural to ask whether vanishing orders of the associated
holomorphic differentials similarly constrain the geometry
of maximal surfaces and their asymptotic boundaries in
higher codimension.

\begin{samepage}
\paragraph{Organization}
Section~\ref{sec:prelim} recalls the affine-sphere and CMC realizations
and establishes the nodal counting lemma.
Section~\ref{sec:pick} proves the projective disk comparison using
regular polygonal models and conormal duality, and derives the
Hilbert--Blaschke norm estimates.
Section~\ref{sec:hopf} develops the support-function comparison for
CMC graphs and proves the hyperbolic disk bound and its rigidity.
Section~\ref{sec:applications} applies a common convex-geometric
estimate to obtain quantitative fence bounds for polygons and
general convex domains, including the critical and equality cases.
Finally, Section~\ref{sec:examples} constructs examples with simple
zeros illustrating the sufficient nature of the multiplicity conditions.
\end{samepage}

\subsection*{Declaration of AI use}
After the second author's preprint \cite{Wang} appeared,
the first author shared a previously obtained, independent
AI-assisted proof of the Fence conjecture and an analogous
result for Hopf zeros. The authors then combined their work and decided to write the paper together.
Comparing the two approaches led them to introduce higher fences
and incorporate duality into the proofs.
Repeated discussions with ChatGPT 6 Astra contributed to the development
of sharp estimates for higher-order zeros and the resulting
higher-fence bounds. Starting from an initial draft produced by ChatGPT,
the authors checked and revised every argument, simplified
the proofs, and reorganized the manuscript.
The authors take full responsibility for the manuscript. 

\subsection*{Acknowledgement}
Both authors are grateful to the Tsinghua Sanya International Mathematics Forum for hosting the conference \emph{Surface Theory: Geometry, Analysis, and Experiments} (25--29 January 2026), which brought together researchers in Higgs bundles, nonabelian Hodge
theory, differential geometry, and integrable systems. Q. Li is partially supported by the National Key R\&D Program of China No. 2022YFA1006600, the Fundamental Research Funds for the Central Universities, and Nankai Zhide foundation. E. Wang is supported by the National Natural Science Foundation of China NSFC--RGC grant No.  12361161603. 

% Add a Declaration of AI here before submission, if required.
% Add acknowledgements here before submission.

\section{Preliminaries}\label{sec:prelim}
We fix the conventions for the two geometric realizations and record the comparison tools used below. We use $\Delta=\partial_x^2+\partial_y^2=4\partial_z\partial_{\bar z}$.

\subsection{Hyperbolic affine spheres and cubic differentials}\label{subsec:affine-prelim}
Let $F:X\to\R^3$ be a hyperbolic affine sphere centered at $0$, normalized so that its Blaschke normal is $F$. Its affine connection $\nabla$ and positive definite \emph{Blaschke metric} $h$ are given by
\begin{equation}\label{eq:affine-gauss}
 D_V(F_*W)=F_*(\nabla_VW)+h(V,W)F,
\end{equation}
where $D$ is the flat ambient connection. This is the affine mean-curvature normalization $H_{\mathrm{aff}}=-1$ of \cite[\S4.1]{DW}.

Let $\nabla^h$ be the Levi-Civita connection of $h$. The difference of the two connections defines the \emph{Pick tensor}
\[
 A(V,W,Z)=h\bigl((\nabla-\nabla^h)_VW,Z\bigr).
\]
We use the convention $A=\Rea C$, where $C$ is the holomorphic Pick differential on the Riemann surface determined by $h$.

On a simply connected surface, compatible data $(h,C)$ determine $F$ up to $\mathrm{SL}(3,\R)$; see \cite[\S4.2 and Proposition~4.6]{DW}. If $h$ is complete, radial projection identifies the properly embedded sphere with a properly convex domain $\Omega\subset\RP$. Conversely, each such domain admits a unique normalized complete affine sphere by the Cheng--Yau theorem \cite{CY,DW}.

In a bounded affine chart, the same sphere has the support-function description
\begin{equation}\label{eq:affine-pde}
 F(x)=-\frac{1}{U(x)}(x_1,x_2,1),\qquad
 \det D^2U=U^{-4},\quad U<0,\quad D^2U>0,\quad U|_{\partial\Omega}=0.
\end{equation}
Here $U$ is smooth in $\Omega$ and continuous on $\overline\Omega$ \cite[\S4.3]{DW}. Differentiating the parametrization gives
\begin{equation}\label{eq:affine-connections}
 h_{ij}=-\frac{U_{ij}}U,
 \qquad
 \Gamma^k_{ij}=-\frac{U_i}{U}\delta^k_j-\frac{U_j}{U}\delta^k_i.
\end{equation}

\subsection{Spacelike CMC surfaces and quadratic differentials}\label{subsec:CMC-prelim}
Equip Minkowski space $\R^{2,1}$ with the metric
$\langle\cdot,\cdot\rangle_L=dx_1^2+dx_2^2-dt^2$.
Increasing $t$ is the future time orientation.

An immersion $\sigma:X\to\R^{2,1}$ is \emph{spacelike} when its
induced metric $I$ is positive definite. Its \emph{future unit normal}
is the timelike unit normal $N$ with positive $t$-component, so it
takes values in the hyperboloid
\[
 \HH=\{\nu\in\R^{2,1}:\langle\nu,\nu\rangle_L=-1,\ \nu_3>0\},
\]
with its curvature-$-1$ metric. Define the shape operator by $\sigma_*B(V)=D_VN$, the second fundamental form by $\II(V,W)=I(BV,W)$, and the mean curvature by $H_{\mathrm{CMC}}=\tfrac12\tr B$. We use strictly convex surfaces with $B>0$ and normalize $H_{\mathrm{CMC}}=1$.

The Gauss map $N:(X,[I])\to\HH$ is harmonic
\cite[Proposition~2.1(vi)]{HTTW}. In our normalization,
\begin{equation}\label{eq:Gauss-Hopf-prelim}
 N^*g_{\HH}=I(B\cdot,B\cdot)=2\II-(\det B)I,
 \qquad Q_N=2(\II)^{2,0}.
\end{equation}
The middle equality follows from $B^2-2B+(\det B)\Id=0$; in particular, $Q_N=0$ exactly at umbilics $B=\Id$.
The Gauss equation and $B>0$, $\tr B=2$ give
\begin{equation}\label{eq:CMC-curvature}
 K_I=-\det B\in[-1,0).
\end{equation}

\paragraph{The CMC realization.}
Let $f:X\to\Omega$ satisfy the hypotheses of Theorem~\ref{thm:hopf}.
Since $\Omega$ is convex and $f$ is a diffeomorphism, $X$ is simply
connected. Write $Q=Q_f=(f^*g_{\HH})^{2,0}$ and define the
$\partial$-energy metric in local coordinates by
\[
 I_f=\rho(f)^2|f_z|^2|dz|^2,
 \qquad g_{\HH}=\rho(w)^2|dw|^2.
\]
The CMC realization of \cite[p.~89 and \S2]{HTTW} gives a
strictly convex spacelike CMC-$1$ immersion, unique up to an ambient
Lorentz isometry and translation, with
\begin{equation}\label{eq:CMC-data}
 I=I_f,\qquad \II=I+\Rea Q,
 \qquad B=I^{-1}\II.
\end{equation}
Here $\Rea(q\,dz^2)=\tfrac12(q\,dz^2+\overline q\,d\overline z^2)$.
Choose the realization so that its Gauss map is $N=f$; then $Q_N=Q$.
Completeness of $I=I_f$ makes the immersion an entire graph
\cite[Proposition~1.2]{BSS}:
\begin{equation}\label{eq:CMC-graph}
 \Sigma_u=\{(x,u(x)):x\in\R^2\},\qquad
 |Du|<1,\quad D^2u>0,\quad
 \operatorname{div}\frac{Du}{\sqrt{1-|Du|^2}}=2.
\end{equation}
Parametrize the graph by $\sigma(x)=(x,u(x))$ and put $W=\sqrt{1-|Du|^2}$.
Its future unit normal and fundamental forms are
\[
 N=\frac{(Du,1)}W,\qquad
 g=(I_{ij})=\Id-Du\otimes Du,\qquad
 \II=W^{-1}D^2u.
\]
The sign follows from $\II_{ij}=-\langle N,\sigma_{ij}\rangle_L=u_{ij}/W$.
Thus $B=g^{-1}\II$ gives
\begin{equation}\label{eq:graph-shape}
 D^2u=WgB
 =\sqrt{1-|Du|^2}\,(\Id-Du\otimes Du)B.
\end{equation}
The Klein coordinate of $N$ is $Du$ \cite[\S2, p.~96]{HTTW}.
Strict convexity makes $Du$ injective. Its image is convex by
Choi--Treibergs \cite[Theorem~4.8]{ChoiTreibergs}; see also
\cite[Proposition~2.1(v)]{HTTW}. Thus $N$ is a harmonic
diffeomorphism onto its image.

For the original harmonic map $f$, completeness of $I=I_f$ and
$N=f$ therefore give
\begin{equation}\label{eq:gradient-image}
 Du(\R^2)=\Omega.
\end{equation}

\subsection{Nodal domains and cyclic symmetry}\label{sec:nodal}
The following counting lemma is the key analytic tool. 

\begin{lemma}\label{lem:nodal}
Let $D\subset\R^2$ be a convex open set and let $w\in C^\infty(D)$ be a nonzero solution of
\begin{equation}\label{eq:linear}
 a^{ij}w_{ij}+b^iw_i-cw=0,
\end{equation}
where the coefficients are smooth, $(a^{ij})$ is locally uniformly positive definite, and $c\ge0$. If $w$ vanishes to order at least $k\ge2$ at an interior point, then $\{w<0\}$ has at least $k$ connected components. If, in addition, $D$ is bounded and $w\in C^0(\overline D)$, then $\{w<0\}\cap\partial D$ has at least $k$ connected components.
\end{lemma}

\begin{proof}
On a small disk choose a positive solution $\phi$ of the same equation, with boundary value $1$. The equation for $w/\phi$ has no zeroth-order term. The Hartman--Wintner local theorem \cite{HartmanWintner}, in the form recalled in \cite[p.~231]{Alessandrini}, gives finite vanishing order $\ell\ge k$ and, after a linear coordinate change, the leading term
\[
 r^\ell\bigl(A\cos\ell\vartheta+B\sin\ell\vartheta\bigr),
 \qquad A^2+B^2>0.
\]
The accompanying first-derivative expansion gives $\ell$ positive sectors alternating with $\ell$ negative sectors. The same local theorem gives unique continuation, so a globally nonzero solution cannot vanish on a neighborhood of the point. Division by $\phi$ does not change the order or signs.

Suppose two different local negative sectors belonged to one negative component. Join small radial arcs in these sectors by a path in that component. By deleting loops and taking the last and first intersections with the radial arcs, one obtains a simple closed curve through the zero on which $w\le0$. Either angular interval between the two negative sectors contains a positive sector. The bounded component enclosed by the curve therefore contains a point where $w>0$. It is compactly contained in $D$: the convex hull of the compact curve is contained in the convex open set $D$. The maximum principle contradicts the existence of a positive interior maximum. Thus the $\ell$ negative sectors belong to distinct components.

Now assume $D$ bounded. On the compact closure of each negative component, $w$ attains a negative minimum. It cannot attain that minimum inside $D$, by the strong maximum principle, and hence reaches a boundary point where $w<0$. At such a point, continuity and the local convexity of $D$ give a connected negative neighborhood inside $D$. The component incident to the boundary is locally constant along the negative boundary set, so each connected component of that set is incident to at most one negative nodal domain. This proves the boundary count.
\end{proof}

We also record a symmetry observation used for the regular polygonal
models.

\begin{lemma}\label{lem:rotation}
Let a conformal automorphism of finite order $k$ fix a point and have primitive $k$-th root of unity as derivative there. An invariant holomorphic $r$-differential, with $k\ge r$, has vanishing order at least $k-r$ at that point, unless it vanishes identically.
\end{lemma}

\begin{proof}
A local holomorphic coordinate obtained by averaging linearizes the automorphism to $z\mapsto\zeta z$. Invariance of $\sum_{j\ge0}a_jz^jdz^r$ gives $a_j\ne0$ only if $k$ divides $j+r$. The least possible nonnegative $j$ is $k-r$.
\end{proof}
\section{Pick differentials and projective disk comparison}\label{sec:pick}
\subsection{Normalization}
We first normalize the domain $\Omega$ at the point $p$. An affine chart is \emph{normalized at $p$} if $p=0$ and its canonical potential satisfies $U(0)=-1$, $DU(0)=0$, and $D^2U(0)=\Id$.

\begin{lemma}\label{lem:affine-contact}
At any prescribed point $p\in\Omega$, a projective change of affine chart, induced by an element of $\mathrm{SL}(3,\R)$, gives an affine chart normalized at $p$ in which $\Omega$ is bounded.
In these coordinates $h_0=|dx|^2$, and the normalized chart is unique
up to an orthogonal linear change of coordinates.
If the Pick differential $C$ has vanishing order at least $m\ge0$ at $p$, then in any such normalized chart the canonical potential $U$ satisfies
\begin{equation}\label{eq:affine-contact}
 U=-\sqrt{1-|x|^2}+O(|x|^{m+3}).
\end{equation}
\end{lemma}

\begin{proof}
Start in a bounded affine chart, translate $p$ to $0$, and write
$U(0)=-s$, $DU(0)=g$, and $D^2U(0)=G$. Then $s>0$ and $G>0$, while
$\det D^2U=U^{-4}$ gives $\det G=s^{-4}$.

Put $\ell(x)=s-g\cdot x$. Convexity gives $U\ge-\ell$, hence
$\ell\ge-U>0$ in $\Omega$. For $a\in\partial\Omega$, the function
\[
 f_a(t)=U(ta)+s-tg\cdot a
\]
satisfies $f_a(0)=f_a'(0)=0$ and $f_a''(t)>0$ for $0<t<1$.
Since $U(a)=0$, it follows that $\ell(a)=\lim_{t\to1}f_a(t)>0$.
Thus $\ell$ has a positive lower bound on the compact set $\overline\Omega$.

Let $B=(sG)^{1/2}$ be the positive definite symmetric square root.
Since $\det B=s^{-1}$, the matrix
\[
 M=\begin{pmatrix}B&0\\-g^T&s\end{pmatrix}
 \in\mathrm{SL}(3,\R)
\]
acts, by the parametrization in \eqref{eq:affine-pde}, as
\[
 y=\frac{Bx}{\ell(x)},\qquad
 \widetilde U(y)=\frac{U(x)}{\ell(x)}.
\]
The lower bound for $\ell$ keeps the new domain bounded.
Since $U+\ell$ has zero value and gradient at $0$ and $D_xy(0)=s^{-1}B$,
differentiation gives
\[
 \widetilde U(0)=-1,\qquad D_y\widetilde U(0)=0,\qquad
 D_y^2\widetilde U(0)=sB^{-T}GB^{-1}=\Id.
\]
Equation~\eqref{eq:affine-connections} then gives $h_0=|dy|^2$.

For uniqueness, a transition between two normalized charts fixes
$0$ and has the form $y=Ax/(c\cdot x+d)$. Choose a homogeneous
representative of absolute determinant $1$ with positive denominator
on the domain. The potential
transforms as $\widetilde U(y)=U(x)/(c\cdot x+d)$.
The conditions $\widetilde U(0)=U(0)=-1$ and
$D\widetilde U(0)=DU(0)=0$ give $d=1$ and $c=0$.
The Hessian conditions then give $A^{-T}A^{-1}=\Id$, so $A$ is
orthogonal.

For the contact estimate, work in any normalized affine chart, writing $x$ for its coordinates and $U$ for its canonical potential.
Use the convention $A=h(\nabla-\nabla^h)$ for the real Pick tensor, so $A=\Rea C$ and $\nabla h=-2A$ \cite[\S4.1]{DW}. Equation~\eqref{eq:affine-connections} yields
\[
 (\nabla_k h)_{ij}
 =-\frac{U_{ijk}}U
  -\frac{U_iU_{jk}+U_jU_{ik}+U_kU_{ij}}{U^2}.
\]
Thus
\begin{equation}\label{eq:pick-tensor}
 A_{ijk}=\frac{1}{4U^2}\partial_i\partial_j\partial_k(U^2).
\end{equation}
Since $C$ has vanishing order at least $m$, its real tensor is $O(|x|^m)$ in these smooth coordinates. Hence $D^3(U^2)=O(|x|^m)$. The normalization conditions and Taylor's theorem give
\[
 U^2=1-|x|^2+O(|x|^{m+3}).
\]
Taking the negative square root gives \eqref{eq:affine-contact}. For $m=0$, no vanishing assumption on $C$ is required; the expansion follows directly from the normalization conditions.
\end{proof}

\subsection{The regular model and its normalization}\label{subsec:affine-model}\label{subsec:polarity}
For $k\ge3$, scale a regular $k$-gon centered at $0$ so that its canonical potential $V$ satisfies $V(0)=-1$, and denote the resulting polygon by $Q_k$. The scaling rule
\[
 V_{\lambda Q}(x)=\lambda^{2/3}V_Q(x/\lambda),\qquad \lambda>0,
\]
makes this scale unique. Rotational symmetry gives $DV(0)=0$ and a scalar Hessian; the equation $\det D^2V=V^{-4}$ then gives $D^2V(0)=\Id$. Thus the chart is normalized at $0$.

To determine its circumradius $R_k$, recall the conormal description of the canonical potential on the polar
\[
 Q_k^\circ=\big\{y\in\R^2:\sup_{x\in\overline Q_k}x\cdot y<1\big\}.
\]
Writing $\ell(x)=x\cdot DV(x)-V(x)\ge1$, the conormal identities are
\[
 y=\frac{DV(x)}{\ell(x)},\qquad V^\circ(y)=-\frac1{\ell(x)}.
\]
At $0$ these formulas give $y(0)=0$, $Dy(0)=\Id$, and the same normalized jet for $V^\circ$. The polar is a rotated regular $k$-gon of circumradius $1/(R_k\cos(\pi/k))$. Uniqueness of the normalized scale therefore gives
\[
 R_k=\frac1{R_k\cos(\pi/k)}.
\]
Hence its circumradius and inradius are
\begin{equation}\label{eq:Rk}
 R_k=\frac1{\sqrt{\cos(\pi/k)}},\qquad
 r_k=R_k\cos(\pi/k)=\sqrt{\cos(\pi/k)}.
\end{equation}

By \cite[Theorem~A]{DW}, the Pick differential is a polynomial cubic differential of degree $k-3$. Its invariance under rotation of order $k$ and Lemma~\ref{lem:rotation} give vanishing order at least $k-3$ at the center, hence exactly $k-3$. This includes the triangle, whose Pick differential is a nonzero constant. For every rotation $Q_{k,\alpha}$ of $Q_k$, its potential $V_\alpha$ therefore satisfies, by Lemma~\ref{lem:affine-contact},
\begin{equation}\label{eq:affine-model-contact}
 V_\alpha=-\sqrt{1-|x|^2}+O(|x|^k).
\end{equation}

\subsection{Proof of Theorem~\ref{thm:pick}}
Normalize at $p$, put $k=m+3$, and fix $\alpha$. On the bounded convex intersection $D=\Omega\cap Q_{k,\alpha}$ set $w=U-V_\alpha$. By Lemma~\ref{lem:affine-contact} and \eqref{eq:affine-model-contact}, $w=O(|x|^k)$. Linearizing the equivalent equation $\log\det D^2U+4\log(-U)=0$ along $U_t=V_\alpha+tw$ gives
\begin{equation}\label{eq:affine-linear}
 a^{ij}w_{ij}-cw=0,\qquad
 a^{ij}=\int_0^1\bigl((D^2U_t)^{-1}\bigr)^{ij}\,dt,
 \quad c=4\int_0^1\frac{dt}{-U_t}>0.
\end{equation}
These coefficients are smooth and locally elliptic because $D^2U_t>0$ and $U_t<0$.

If $w\equiv0$, the boundary values and interior negativity of the two
potentials force $\Omega=Q_{k,\alpha}$; the regular-case radii are then
given by \eqref{eq:Rk}. We may therefore assume $w\not\equiv0$ for
every $\alpha$.

The negative boundary set is exactly
\begin{equation}\label{eq:boundary-sign}
 \{w<0\}\cap\partial D=\partial Q_{k,\alpha}\cap\Omega.
\end{equation}
Lemma~\ref{lem:nodal} gives at least $k$ components. Each of the $k$ closed sides of $Q_{k,\alpha}$ meets the convex set $\Omega$ in at most one interval, so there are at most $k$. Hence every side meets the open domain $\Omega$.

Let $h_\Omega(\nu)=\max_{x\in\overline\Omega}\nu\cdot x$ be the
support function of $\Omega$. For each $\nu\in S^1$, rotate the model
so that one side lies on $\nu\cdot x=r_k$.
Since this side meets the open set $\Omega$, we have
$h_\Omega(\nu)>r_k$.
Continuity and compactness of $S^1$ give
$\min_{S^1}h_\Omega>r_k$, proving $B_{r_k}\Subset\Omega$.

For the outer inclusion, consider the polar $\Omega^\circ$. The conormal formulas in Section~\ref{subsec:polarity}, with $V$ replaced by $U$, give the same normalized jet at $0$. Conormal duality preserves the Blaschke metric and changes $C$ to $-C$ \cite[\S2.1.6]{LM}, so its vanishing order is unchanged. Since polarity exchanges the normalized regular $k$-gons, $\Omega^\circ$ is also nonregular. The strict inner inclusion just proved gives $ B_{r_k}\Subset\Omega^\circ$. By polarity,
\[
 \Omega\Subset B_{1/r_k}.
\]
 The regular models establish sharpness.
\hfill\qed

\subsection{The Hilbert norm and vanishing order}
Benoist--Hulin proved uniform comparability of the Hilbert and
Blaschke metrics \cite{BHFinite}; Hildebrand's estimates give the
constants $1/\sqrt2$ and $\sqrt2$ below at $m=0$
\cite[Corollary~1.5]{Hildebrand}. Vanishing order refines these
constants, with multiplicative error $1+O(m^{-2})$ as $m\to\infty$.

We recall the Hilbert metric. If $x,y\in\Omega$ are collinear with boundary points $a,x,y,b$ in that order, then
\[
 d_\Omega(x,y)=\frac12\log\frac{|y-a|\,|x-b|}{|x-a|\,|y-b|}.
\]
At $0$, if the two endpoints along a unit direction have Euclidean distances $\rho_+$ and $\rho_-$, its infinitesimal norm is
\begin{equation}\label{eqn:HilbertNorm}
 F_\Omega(0,v)=\frac{|v|}{2}\left(\frac1{\rho_+}+\frac1{\rho_-}\right).
\end{equation}

\begin{corollary}\label{cor:Hilbert-norm}
At a point $p$ where the Pick differential has finite vanishing order $m\ge0$, put $k=m+3$.
Then, for every $v\in T_p\Omega$,
\begin{equation}\label{eq:Hilbert-norm}
 \sqrt{\cos\frac{\pi}{k}}\,\|v\|_h
 \le F_\Omega(p,v)
 \le
 \frac{\|v\|_h}{\sqrt{\cos(\pi/k)}}.
\end{equation}
When $k$ is even, both norm constants are sharp, as realized at the center of the normalized regular $k$-gon.
\end{corollary}
\begin{proof}
Use normalized coordinates at $p$, so $p=0$ and $\|v\|_h=|v|$. For $v\ne0$, Theorem~\ref{thm:pick} gives $r_k\le\rho_\pm\le R_k$. Substitution into \eqref{eqn:HilbertNorm} yields $|v|/R_k\le F_\Omega(0,v)\le|v|/r_k$, proving \eqref{eq:Hilbert-norm}. In the regular even $k$-gon, directions through opposite vertices have $\rho_+=\rho_-=R_k$, while directions through opposite side midpoints have $\rho_+=\rho_-=r_k$; these attain the lower and upper bounds, respectively.
\end{proof}

\section{Hopf differentials and hyperbolic disk comparison}\label{sec:hopf}
\subsection{Normalization}\label{subsec:hopf-normalization}
Let $f$ satisfy Theorem~\ref{thm:hopf}. Section~\ref{subsec:CMC-prelim} realizes it as the Gauss map of the complete strictly convex graph \eqref{eq:CMC-graph}, with gradient image \eqref{eq:gradient-image}. We call a graph \emph{normalized at a specified point} if this point lies over $x=0$ and $u(0)=0$, $Du(0)=0$.

\begin{lemma}\label{lem:CMC-contact}
Let $\Sigma$ be a complete strictly convex spacelike CMC-$1$ surface, and let $Q_N$ be the Hopf differential of its Gauss map $N$. At any prescribed point, a time-orientation-preserving Lorentz isometry and an ambient translation give an entire graph $t=u(x)$ normalized at that point, whose Gauss map takes the point to the center of the Klein disk. If $Q_N$ has vanishing order at least an integer $m\ge0$ there, then in any such normalization
\begin{equation}\label{eq:CMC-contact}
 u(x)=\sqrt{1+|x|^2}-1+O(|x|^{m+2}).
\end{equation}
\end{lemma}

\begin{proof}
Choose a time-orientation-preserving Lorentz isometry sending the future unit normal at the prescribed point to $(0,0,1)$, then translate that point to the origin. The transformed surface is complete and hence an entire graph by \cite[Proposition~1.2]{BSS}. Its tangent plane at the origin is horizontal, so $u(0)=0$ and $Du(0)=0$.

For $m=0$, \eqref{eq:CMC-contact} follows directly from these two conditions and Taylor's theorem, without a vanishing assumption on $Q_N$. Assume now that $m\ge1$. In any normalized graph coordinates, \eqref{eq:Gauss-Hopf-prelim} and $\tr B=2$ give $\II-I=\Rea Q_N$. Thus the assumed vanishing order implies
\[
 B-\Id=I^{-1}\Rea Q_N=O(|x|^m),
\]
since the induced metric is smooth and positive definite. By \eqref{eq:graph-shape},
\[
 D^2u=A(Du)+O(|x|^m),\qquad
 A(p)=\sqrt{1-|p|^2}\,(\Id-p\otimes p),
\]
where $p\in\R^2$, $|p|<1$, is the gradient variable, evaluated at
$p=Du(x)$ in $A(Du)$. Multiplication by the smooth bounded matrix
$A(Du)$ preserves the order of the error.
The hyperboloid $u_0=\sqrt{1+|x|^2}-1$ satisfies $D^2u_0=A(Du_0)$ with the same value and gradient at $0$. Put $w=u-u_0$. If its Taylor polynomial through degree $m+1$ were nonzero, let $P_\ell$ be its lowest-degree homogeneous term, where $2\le\ell\le m+1$. Smoothness of $A$ gives
\[
 D^2w=A(Du)-A(Du_0)+O(|x|^m)
      =O(|Dw|)+O(|x|^m)=O(|x|^{\ell-1}).
\]
But $D^2w=D^2P_\ell+O(|x|^{\ell-1})$, so the homogeneous polynomial $D^2P_\ell$ of degree $\ell-2$ must vanish. This contradicts $P_\ell\ne0$. Taylor's theorem therefore gives $w=O(|x|^{m+2})$.
\end{proof}

\subsection{Regular polygonal models}

For $k\ge3$, let $b_1,\ldots,b_k\in S^1$ be equally spaced.
By \cite[Propositions~2.1--2.2]{HTTW}, there is a unique entire,
complete, strictly convex spacelike CMC-$1$ graph whose support
function is zero at $b_1,\ldots,b_k$ and $+\infty$ elsewhere
on $S^1$. After a vertical translation, denote its graph function
by $v$, normalized by $v(0)=0$.
By \eqref{eq:support-domain}, its Gauss image is the regular ideal
$k$-gon with vertices $b_1,\ldots,b_k$.
Uniqueness gives invariance under rotation by $2\pi/k$, hence
$Dv(0)=0$.
The rotational symmetry and Lemma~\ref{lem:rotation} give Hopf
vanishing of order at least $k-2$ at the center.
Lemma~\ref{lem:CMC-contact} therefore yields
\[
    v(x)=\sqrt{1+|x|^2}-1+O(|x|^k).
\]

\subsection{Support functions and Legendre duality}\label{subsec:legendre-support}
For an entire strictly convex spacelike CMC-$1$ graph
$\Sigma_w=\operatorname{graph}w$, its support function in Klein
coordinates is the Legendre transform
\[
    w^*(p)=\sup_{q\in\Sigma_w}\langle q,(p,1)\rangle_L
          =\sup_{x\in\mathbb R^2}\bigl(p\cdot x-w(x)\bigr),
    \qquad |p|\le1.
\]
It is convex and lower semicontinuous, with values in
$\mathbb R\cup\{+\infty\}$.

On the open Gauss image $Dw(\mathbb R^2)$, the Legendre transform can also be written
using the gradient coordinates $p=Dw(x)$. Indeed, convexity gives
$w(z)\ge w(x)+Dw(x)\cdot(z-x)$, so the supremum is attained at
$x$ when $p=Dw(x)$. Thus
\[
    w^*(Dw(x))=x\cdot Dw(x)-w(x).
\]
The normal-translation correspondence of \cite[Corollary~5.5]{BSS},
expressed in our normalization $H_{\mathrm{CMC}}=1$, shows that
the past normal translate $\sigma-\tfrac12N$ of $\Sigma_w$
is an entire surface of intrinsic Gaussian curvature $-4$,
with the same Gauss image. Denote its support function by $s_w$.
Evaluating the normal displacement gives
\begin{equation}\label{eq:normal-support}
    s_w(p)=w^*(p)+\frac12\sqrt{1-|p|^2},
    \qquad |p|\le1.
\end{equation}
 In particular, $s_w$ and $w^*$ have the same
finite domain and agree on $S^1$.
The function $s_w$ is convex and lower
semicontinuous, and is smooth with positive definite Hessian on
$Dw(\mathbb R^2)$. Since the translated surface has
$\det B=4$, \cite[Lemma~2.11]{BSS19} gives
\begin{equation}\label{eq:support-MA}
    \det D^2s_w(p)=\frac{1}{4(1-|p|^2)^2},
    \qquad p\in Dw(\mathbb R^2).
\end{equation}
Applying \cite[Theorem~4.4]{BSS19} to $s_w$, we obtain
\begin{equation}\label{eq:support-domain}
    \{w^*<+\infty\}
       =\operatorname{conv}\{b\in S^1:w^*(b)<+\infty\},
    \qquad
    Dw(\mathbb R^2)=\operatorname{int}\{w^*<+\infty\}.
\end{equation}
Moreover, on any boundary chord $[b_-,b_+]$ of $Dw(\mathbb R^2)$,
with $b_\pm\in S^1$, the same theorem gives the following
alternative: if both endpoint values of $s_w$ are finite, then
$s_w$ is their affine interpolation on the chord; if either
endpoint value is infinite, then $s_w=+\infty$ on the open chord.

We will also use the elementary fact that a finite lower
semicontinuous convex function on a compact polygon is continuous.
Indeed, a triangulation having a prescribed point as a vertex
gives a continuous piecewise affine upper bound agreeing with
the function at that point; lower semicontinuity gives the
opposite inequality in the limit.

\subsection{Counting comparison components}
We use the support functions to bound the number of comparison
components from above. Although the model's dual domain is a compact
polygon, the other support function may take the value $+\infty$
there. Thus the boundary version of Lemma~\ref{lem:nodal} does not
apply directly; the convex patching argument below handles these
possibly infinite boundary values.

\begin{lemma}\label{lem:comparison-components}
Let $\Sigma_u,\Sigma_v$ be complete strictly convex CMC-$1$ graphs.
Suppose that
\[
 Dv(\mathbb R^2)=\operatorname{int}P,\qquad
 P=\operatorname{conv}\{b_1,\ldots,b_k\},
\]
where $b_1,\ldots,b_k\in S^1$ are distinct and $k\ge3$.
Then $\{u<v\}$ has at most $k$ connected components.
Moreover, if $Du(\mathbb R^2)$ is disjoint from some edge of $P$,
then $\{u<v\}$ has at most $k-1$ components.
\end{lemma}

\begin{proof}
By \eqref{eq:support-domain}, the points of $S^1$ where $v^*$
is finite lie in $P\cap S^1=\{b_1,\ldots,b_k\}$, and their
convex hull has interior $\operatorname{int}P$.
Hence all vertices belong to this set, and the finite domain
of $v^*$ is exactly $P$.
The function $s_v$ is finite and continuous on $P$, and affine
on each edge. By \eqref{eq:normal-support}, $v^*$ is also
continuous on $P$, and
\[
    A:=\{p\in P:u^*(p)>v^*(p)\}
      =\{p\in P:s_u(p)>s_v(p)\}.
\]
For $p\in P$, put
\[
    C_p=\{x\in\mathbb R^2:u(x)<p\cdot x-v^*(p)\}.
\]
Each $C_p$ is convex and contained in $\{u<v\}$, and it is
nonempty precisely when $p\in A$.
We record a connectivity observation used again below:
for every connected $K\subset A$, the set
\begin{equation}\label{eq:dual-connected-union}
    \bigcup_{p\in K} C_p
\end{equation}
is connected and contained in $\{u<v\}$.
Indeed, a point of $C_p$ remains in $C_q$ for all nearby
parameters $q$, by continuity of $v^*$. Since each $C_p$ is
connected, the component of the union containing $C_p$ is
locally constant in $p\in K$, and hence constant.

We claim that every connected component of $A$ contains a
vertex of $P$. Suppose that a component $A_0$ contains no
vertex, and define
\[
    h(p)=
    \begin{cases}
        s_u(p),&p\in A_0,\\
        s_v(p),&p\in P\setminus A_0.
    \end{cases}
\]
The set $A$ is relatively open in $P$, as are its components.
The function $h$ is lower semicontinuous: on $A_0$ this follows
from $s_u$, and elsewhere from $h\ge s_v$ and continuity of $s_v$.
It is also convex. To see this, consider a segment whose
endpoint values of $h$ are finite, and let $\ell$ be their
affine interpolation. Convexity gives $s_v\le\ell$ along the
segment. Consider an interval of its intersection with $A_0$.
At an interval endpoint outside $A_0$, we have
$s_u\le s_v\le\ell$, since a boundary point of $A_0$ cannot
belong to $A$. An interval endpoint in $A_0$ must be an endpoint
of the original segment, where $s_u=h=\ell$.
Convexity of $s_u$ therefore gives $s_u\le\ell$ throughout the
interval. Thus $h\le\ell$ on the whole segment, proving
convexity even when $s_u$ takes infinite values.

At every vertex $h=s_v<+\infty$. Convexity therefore makes
$h$ finite throughout $P$, and the elementary fact recorded in
Section~\ref{subsec:legendre-support} makes it continuous on $P$.
Since $s_v$ is affine on each edge, $h\ge s_v$ and equality at
the vertices imply $h=s_v$ on $\partial P$.
For sufficiently small $\varepsilon>0$, the function
\[
    h(p)-s_v(p)+\varepsilon(|p|^2-1)
\]
has a positive maximum at an interior point of $P$.
This point lies in $A_0$, where $h=s_u$ locally.
As $A_0$ is open and $s_u$ is finite there, the point belongs
to $Du(\mathbb R^2)$; both support functions are therefore smooth
at this point. The Hessian inequality at the maximum gives
\[
    D^2s_u+2\varepsilon\Id\le D^2s_v,
\]
contradicting their positive definite Hessians and the common
determinant in \eqref{eq:support-MA}. This proves the claim.

Now let $N$ be a component of $\{u<v\}$ and choose $x\in N$.
For $p=Dv(x)$, Legendre equality gives
\[
    p\cdot x-v^*(p)=v(x)>u(x),
\]
so $x\in C_p$ and $p\in A$.
The component of $A$ containing $p$ contains some vertex $b_j$.
By \eqref{eq:dual-connected-union}, the corresponding union
lies in $N$ and contains $C_{b_j}$.
Thus every component contains at least one nonempty $C_{b_j}$.
Each such set is connected, so there are at most $k$ components.

Finally, suppose that $Du(\mathbb R^2)$ is disjoint from an
edge $E=[b_-,b_+]$ of $P$.
If any $C_{b_j}$ is empty, the preceding argument gives at most
$k-1$ components. Otherwise, nonemptiness of the two endpoint
sets gives
\[
    s_u(b_\pm)=u^*(b_\pm)>v^*(b_\pm)=s_v(b_\pm).
\]
Since $E$ is a full chord of the Klein disk, the convex Gauss
image $Du(\mathbb R^2)$ lies on one side of the line containing
$E$. If $s_u$ is finite at a point of $E^\circ$,
\eqref{eq:support-domain} expresses that point as a convex
combination of finite ideal directions on this line.
These can only be $b_-$ and $b_+$, so both endpoint values
are finite. Then $E$ is a boundary chord of the Gauss image,
and $s_u|_E$ is their affine interpolation.
Otherwise $s_u=+\infty$ throughout $E^\circ$.
As $s_v$ is affine on $E$, the endpoint inequalities give
$E\subset A$ in either case.
The connected union \eqref{eq:dual-connected-union}, with $K=E$,
therefore places $C_{b_-}$ and $C_{b_+}$ in the same component
of $\{u<v\}$. Since every component contains some $C_{b_j}$,
there are at most $k-1$ components.
\end{proof}

\begin{corollary}\label{cor:model-hopf-order}
The regular model $v$ has a Hopf zero of order exactly $k-2$
at the center.
\end{corollary}
\begin{proof}
Its Gauss image is an ideal polygon, so the model is
not totally umbilic and its Hopf differential is not identically zero.
If the central order were at least $k-1$,
Lemma~\ref{lem:CMC-contact} would give
\[
    \sqrt{1+|x|^2}-1-v(x)=O(|x|^{k+1}).
\]
This nonzero difference satisfies the linear elliptic equation
obtained by subtracting the two CMC equations, with no zeroth-order term.
Lemma~\ref{lem:nodal} then gives at least
$k+1$ components of
$\{\sqrt{1+|x|^2}-1<v\}$, whereas
Lemma~\ref{lem:comparison-components}, applied with the hyperboloid
as the first graph, gives at most $k$, a contradiction.
\end{proof}

\subsection{Proof of Theorem C}

\begin{proof}[Proof of Theorem C]
The case $m=0$ follows from $r_0=0$ and $f(z_0)\in\Omega$.
Assume $m\ge1$, and put
\[
    k=m+2,\qquad c=\cos\frac{\pi}{k}.
\]
Use the CMC realization and normalization of
Sections~\ref{subsec:CMC-prelim} and~\ref{subsec:hopf-normalization}:
\[
    u(0)=0,\qquad Du(0)=0,\qquad Du(\mathbb R^2)=\Omega,
    \qquad
    u(x)=\sqrt{1+|x|^2}-1+O(|x|^k).
\]
Klein coordinates are centered at $f(z_0)$.

Suppose that $\overline B_c(0)\not\subset\Omega$.
Separating a point of $\overline B_c(0)\setminus\Omega$ from
the open convex set $\Omega$, choose $\xi\in S^1$ such that
\[
    p\cdot\xi\le c\qquad(p\in\Omega).
\]
Rotate the normalized regular model $v$ so that an edge $E$
of its closed Klein polygon lies on $p\cdot\xi=c$.
Since $\Omega$ is open, it is disjoint from $E$.
If $u\ne v$, the nonzero function $u-v=O(|x|^k)$ satisfies a
smooth locally uniformly elliptic linear equation with no
zeroth-order term. Lemma~\ref{lem:nodal} therefore gives at
least $k$ components of $\{u<v\}$, whereas the last assertion
of Lemma~\ref{lem:comparison-components} gives at most $k-1$.
This contradiction forces $u=v$, so $\Omega$ is the centered
regular ideal $k$-gon.
The normalization can be chosen to preserve the ambient orientation
and time orientation. Since the induced metric of the CMC
realization is $I_f$, the identity $u=v$ identifies the source
conformally with the regular model and identifies $f$ with its
Gauss map, up to an orientation-preserving target isometry.

Hence $\overline B_c(0)\Subset\Omega$, unless $\Omega$ is the
centered regular ideal $k$-gon. Finally,
\[
    d_{\mathbb H^2}(0,p)=\operatorname{arctanh}|p|,
    \qquad
    \operatorname{arctanh}\!\left(\cos\frac{\pi}{k}\right)
       =\log\cot\frac{\pi}{2k}.
\]
The regular ideal $k$-gon has exactly this inradius, proving
the estimate and rigidity.

By Corollary~\ref{cor:model-hopf-order}, the regular model has a
central Hopf zero of order $k-2=m$, so it establishes sharpness
for every $m\ge1$.
For $m=0$, take points in the ideal-triangle model approaching an
interior point of a side, avoiding the images of its isolated
Hopf zeros. Their distance to the boundary tends to zero.
\end{proof}

\section{Fences and higher fences}\label{sec:applications}

We first prove a convex-geometric estimate, then apply it to
higher fences for polygons and general convex domains.

For a compact convex set $K\subset\R^2$ with nonempty interior,
a point $v\in K$ is an \emph{extreme point} if it does not lie
strictly between two distinct points of $K$.
For a polygon these are precisely its vertices; for a strictly
convex set, every boundary point is extreme.

\begin{lemma}\label{lem:chain-disk}
Let $K\subset\R^2$ be a compact convex set with nonempty interior,
and suppose
\[
 B_a(0)\subset\Int K,\qquad K\subset\overline B_b(0),
 \qquad 0<a<b,\qquad \vartheta=\arccos(a/b).
\]
Let $v_0,\ldots,v_j$, with $j\ge2$, be distinct extreme points
occurring in order along a boundary arc of $K$, with each
$[v_{i-1},v_i]\subset\partial K$.
Let $H$ be the closed half-plane bounded by $v_0v_j$ and excluding
$v_1,\ldots,v_{j-1}$.
If $j\vartheta\le\pi/2$, then
\[
 \{x\in\R^2:|x|\le b\cos(j\vartheta)\}
 \subset H\cap\Int K.
\]
\end{lemma}

\begin{proof}
Write $\rho_i=|v_i|$ and
$t_i=\arccos(a/\rho_i)\le\vartheta$.
Each boundary segment lies on a supporting line at distance at
least $a$ from $0$. Its endpoint rays therefore have angular gap
$\delta_i\le t_{i-1}+t_i$.
The total angular span $A$ of the chain satisfies
\[
 A\le 2(j-1)\vartheta+t_0+t_j
 \le2j\vartheta\le\pi.
\]
Consequently, $0\in H$ and \begin{equation}\label{eqn:A}A-(t_0+t_j)\le 2(j-1)\vartheta.\end{equation}

Let $h$ be the distance from $0$ to the line $v_0v_j$.
If $h\ge a$, then
$h\ge a=b\cos\vartheta>b\cos(j\vartheta)$.
Otherwise its perpendicular foot lies in $B_a\subset\Int K$.
Since $v_0,v_j$ are extreme points, the intersection of their line
with $K$ is precisely $[v_0,v_j]$, so the foot lies in the interior
of this segment. Put
$\alpha_i=\arccos(h/\rho_i)$ for $i=0,j$.
Then $A=\alpha_0+\alpha_j$.
The function
\[
 \rho\longmapsto\arccos(h/\rho)-\arccos(a/\rho)
\]
is decreasing on $[a,\infty)$, and hence
\[
 \alpha_i-t_i\ge\arccos(h/b)-\vartheta.
\]
Combining with (\ref{eqn:A}), we have $\arccos(h/b)\le j\vartheta$ and thus
$h\ge b\cos(j\vartheta)$.
Thus the stated disk lies in $H$. Its radius is strictly less
than $a$, so it also lies in $\Int K$.
\end{proof}

\subsection{Fences for polygons}\label{subsec:polygon-fences}
Let $P$ be a convex projective $n$-gon, or an ideal $n$-gon in the
Klein disk, and use the higher fences defined in \eqref{eq:fence-def}.

Although the definition makes sense for $2\le s\le n-1$,
the fence $F_s(P)$ is necessarily empty when $2s>n$.
Thus a nonempty higher fence can occur only for
\[
  2\le s\le \left\lfloor\frac n2\right\rfloor.
\]
At the critical value $n=2s$, the defining diagonals
join opposite vertices, and $F_s(P)$ is either empty or
their single common intersection point.

\begin{corollary}\label{cor:higher-fences}
Let $n\ge3$ and $m\ge0$ be integers. Consider either of the
following settings:
\begin{enumerate}[label=\textup{(\roman*)},leftmargin=2em]
\item $P$ is a convex projective $n$-gon whose Pick differential
has vanishing order $m$ at $p\in\Int P$. Set $k=m+3$, and let
$d=d_P$ be the Hilbert metric with the half-logarithmic
cross-ratio normalization.
\item $f:X\to P$ is an orientation-preserving harmonic
diffeomorphism onto an ideal hyperbolic $n$-gon, with complete
$\partial$-energy metric. Its Hopf differential has vanishing
order $m$ at $z_0\in X$. Set $p=f(z_0)$, $k=m+2$, and
$d=d_{\HH}$.
\end{enumerate}
For every integer $2\le s\le n-1$ with $k\ge2s$, we have
$p\in F_s(P)$. Moreover,
\begin{equation}\label{eq:higher-fence-distance}
 d\bigl(p,\partial F_s(P)\bigr)
 >\log\cot\frac{s\pi}{2k},
\end{equation}
unless one of the following occurs:
\begin{enumerate}[label=\textup{(\alph*)},leftmargin=2em]
\item In the Pick case, $P$ is the regular $k$-gon in the chart
normalized at $p$, and either $k=2s$ or both $k$ and $s$ are even.
\item In the Hopf case, $P$ is a regular ideal $k$-gon centered at $p$.
\end{enumerate}
In each exceptional case, equality holds.
At the critical value $k=2s$, the point $p$ lies in $\Int F_s(P)$
unless $P$ is the corresponding regular $2s$-gon; in that exceptional
case, $F_s(P)=\{p\}$.
\end{corollary}

\begin{proof}
Use the chart normalized at $p$ in the Pick case and the Klein disk
centered at $p$ in the Hopf case. Let $a,b$ be the actual centered
inradius and circumradius, so $b=1$ in the Hopf case, and set
$\vartheta=\arccos(a/b)$.
Theorems~\ref{thm:pick} and~\ref{thm:hopf} give
$a/b\ge\cos(\pi/k)$, hence $\vartheta\le\pi/k$, with strict
inequality outside the corresponding regular $k$-gon.
For $k\ge2s$, Lemma~\ref{lem:chain-disk}, applied to each chain
of $s$ sides in $K=\overline P$, gives
\begin{equation}\label{eq:fence-annulus-disk}
 \{x\in\R^2:|x|\le b\cos(s\vartheta)\}\subset F_s(P).
\end{equation}
In particular, $p\in F_s(P)$.
The defining cuts remove every original vertex and side, so
$F_s(P)\Subset P$ in either setting. Moreover, for $x\in P$,
$d(p,x)\ge\operatorname{arctanh}(|x|/b)$: in the Pick case this
follows from monotonicity of the Hilbert metric under inclusion
$P\subset B_b$, while in the Hopf case it is the hyperbolic radial
formula, with equality. Thus
\begin{equation}\label{eq:fence-annulus-distance}
 d\bigl(p,\partial F_s(P)\bigr)
 \ge\operatorname{arctanh}\cos(s\vartheta)
 =\log\cot\frac{s\vartheta}{2}
 \ge\log\cot\frac{s\pi}{2k}.
\end{equation}
The last inequality is strict when $\vartheta<\pi/k$, since
$\log\cot t$ is strictly decreasing.

\smallskip
\noindent\emph{Equality in the Pick case.}
For equality when $k>2s$, it remains to consider the regular $k$-gon
of circumradius $R=R_k$. Put $\ell=R\cos(s\pi/k)>0$. Its
Euclidean closest fence-boundary points to $0$ are the diagonal midpoints
$x_i=\ell e_i$, where
\[
 \arg e_i=\beta_0+\frac{(2i+s)\pi}{k}
 \quad\text{if}\quad
 \arg v_i=\beta_0+\frac{2i\pi}{k}.
\]
Indeed, every defining line has distance $\ell$ from $0$, and its
foot $x_i$ belongs to the disk in \eqref{eq:fence-annulus-disk}.
Since $\partial F_s(P)$ is compactly contained in $P$, the Hilbert
distance to it is attained. The inequalities
\[
 d_P(0,x)\ge\operatorname{arctanh}(|x|/R)
 \ge\operatorname{arctanh}(\ell/R),\qquad x\in\partial F_s(P),
\]
show that equality in the desired bound requires $x=x_i$ for some $i$
and equality in the disk-metric comparison.
To characterize the latter, write $x=te$ with $t>0$ and $|e|=1$,
and let $-\rho_-e,\rho_+e$ be the endpoints of the chord
$P\cap\R e$. Since $\rho_\pm\le R$, the cross-ratio formula gives
\[
 d_P(0,te)=\frac12\log\frac{1+t/\rho_-}{1-t/\rho_+}
 \ge\frac12\log\frac{1+t/R}{1-t/R}
 =\operatorname{arctanh}(t/R).
\]
Because $t>0$, equality holds exactly when $\rho_- =\rho_+ =R$.
Only vertices of the regular polygon lie on its circumcircle, so both
chord endpoints must be vertices.
At a midpoint $x_i$ this happens exactly when $k$ and $s$
are even: $k$ must be even for opposite vertices to exist, and the
displayed angle is a vertex direction exactly when $s$ is even.
Conversely, under these conditions $x_i$ attains the distance bound.
This proves the stated equality criterion.

\smallskip
\noindent\emph{Equality in the Hopf case.}
In the regular case, the Euclidean midpoints of the defining chords
have radius $\cos(s\pi/k)$. The disk inclusion
\eqref{eq:fence-annulus-disk} places these points on the fence
boundary, so the distance bound is attained.

Finally, at $k=2s$, strictness gives $p\in\Int F_s(P)$ in every
nonregular case. In either regular model, every defining diagonal
is diametral and occurs with both orientations. Their retained
half-planes intersect only at $0$, so $F_s(P)=\{p\}$.
\end{proof}

\subsection{Fences for general convex domains}\label{subsec:general-fences}
Let $\Omega$ be a properly convex projective domain, or a hyperbolic
domain as in Theorem~\ref{thm:hopf}. Let $K$ be its closure in a
bounded affine chart or in the closed Klein disk, respectively.
For an integer $s\ge2$, consider distinct extreme points $v_0,\ldots,v_j$,
$2\le j\le s$, occurring in order along a boundary arc of $K$,
with each $[v_{i-1},v_i]$ a boundary segment. Let
$H(v_0,\ldots,v_j)$ be the closed half-plane bounded by $v_0v_j$
and excluding the intermediate vertices. Define
\[
 F_s(\Omega)=\Omega\cap\bigcap H(v_0,\ldots,v_j),
\]
where the intersection runs over all such chains; if there are
none, set $F_s(\Omega)=\Omega$. This agrees with
\eqref{eq:fence-def} for polygons in its stated range. The extreme-point set and the family of defining chains may be infinite; only the length of each chain is bounded by $s$.

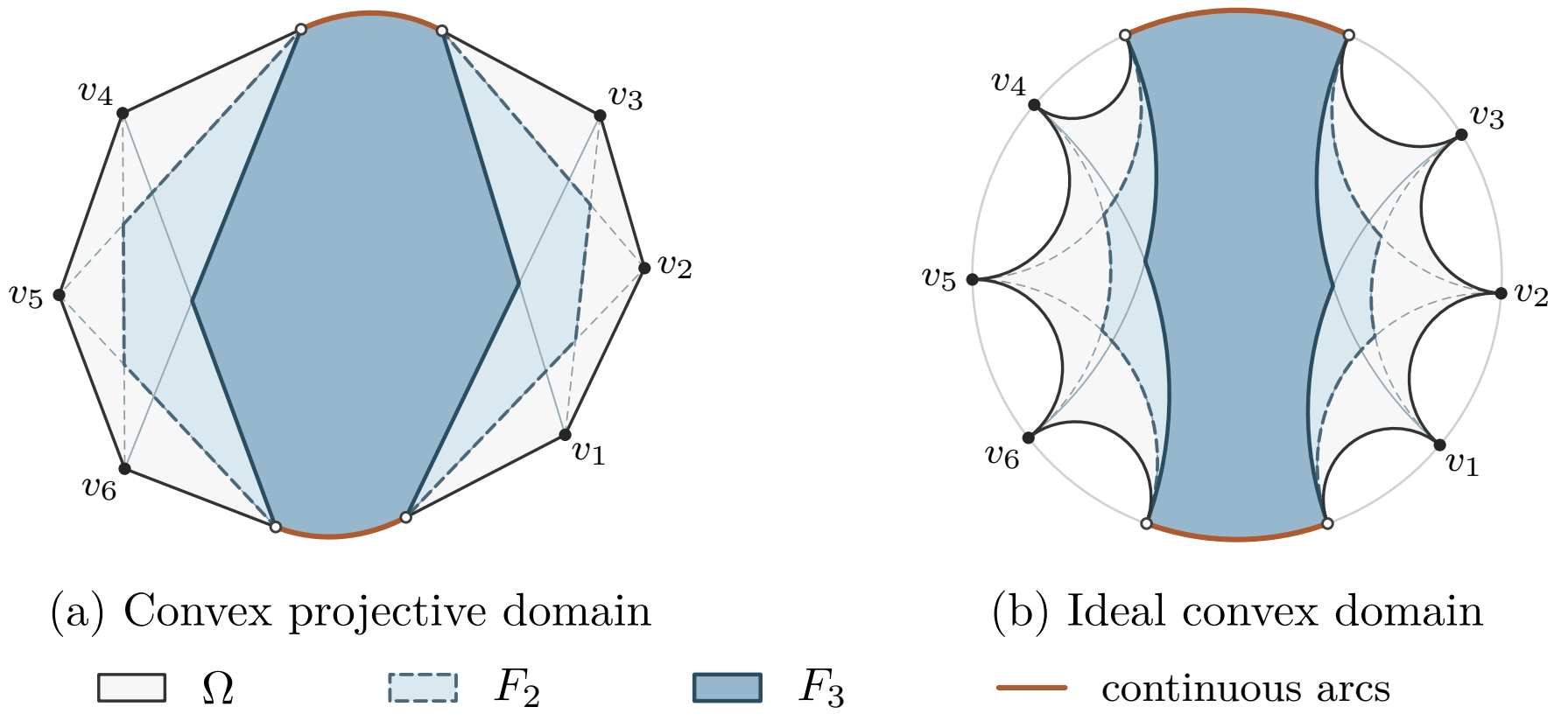
\begin{figure}[H]
\centering
\begin{tikzpicture}[scale=2.0,line cap=round,line join=round,
 every node/.style={font=\small}]
\definecolor{generalFenceFill}{RGB}{218,232,239}
\definecolor{generalFenceDeep}{RGB}{151,184,204}
\definecolor{generalFenceInk}{RGB}{42,77,98}
\definecolor{generalFenceArc}{RGB}{173,91,48}
\path[use as bounding box] (-1.39,-1.68) rectangle (4.72,1.29);
\begin{scope}
\begin{scope}[cm={.75,0,.08,.99,(0,0)}]
\fill[black!3] (-1.211768873,0.617427080) -- (-1.457999121,-0.076410496) -- (-1.056706137,-0.739913603) -- (-0.272659873,-0.962110489) arc[start angle=254.177393686,end angle=292.284862768,radius=1] -- (1.147831871,-0.610313032) -- (1.479774589,0.025829562) -- (1.193948742,0.608347270) -- (0.362733795,0.931892802) arc[start angle=68.731817040,end angle=110.332074506,radius=1] -- (-1.211768873,0.617427080) -- cycle;
\fill[generalFenceFill] (-1.163252762,0.192741612) -- (-1.102126983,-0.342322513) -- (-0.272659873,-0.962110489) arc[start angle=254.177393686,end angle=292.284862768,radius=1] -- (1.161495826,-0.249236506) -- (1.181073558,0.268114378) -- (0.362733795,0.931892802) arc[start angle=68.731817040,end angle=110.332074506,radius=1] -- (-1.163252762,0.192741612) -- cycle;
\fill[generalFenceDeep] (-0.785776471,-0.099072287) -- (-0.272659873,-0.962110489) arc[start angle=254.177393686,end angle=292.284862768,radius=1] -- (0.853595939,-0.032331268) -- (0.362733795,0.931892802) arc[start angle=68.731817040,end angle=110.332074506,radius=1] -- (-0.785776471,-0.099072287) -- cycle;
\draw[generalFenceInk!55,line width=.33pt,dash pattern=on 1.8pt off 1.4pt] (-0.347460633,0.937694571) -- (-1.457999121,-0.076410496);
\draw[generalFenceInk!55,line width=.33pt,dash pattern=on 1.8pt off 1.4pt] (-1.211768873,0.617427080) -- (-1.056706137,-0.739913603);
\draw[generalFenceInk!55,line width=.33pt,dash pattern=on 1.8pt off 1.4pt] (-1.457999121,-0.076410496) -- (-0.272659873,-0.962110489);
\draw[generalFenceInk!55,line width=.33pt,dash pattern=on 1.8pt off 1.4pt] (0.379211711,-0.925309936) -- (1.479774589,0.025829562);
\draw[generalFenceInk!55,line width=.33pt,dash pattern=on 1.8pt off 1.4pt] (1.147831871,-0.610313032) -- (1.193948742,0.608347270);
\draw[generalFenceInk!55,line width=.33pt,dash pattern=on 1.8pt off 1.4pt] (1.479774589,0.025829562) -- (0.362733795,0.931892802);
\draw[generalFenceInk!48,line width=.36pt] (-0.347460633,0.937694571) -- (-1.056706137,-0.739913603);
\draw[generalFenceInk!48,line width=.36pt] (-1.211768873,0.617427080) -- (-0.272659873,-0.962110489);
\draw[generalFenceInk!48,line width=.36pt] (0.379211711,-0.925309936) -- (1.193948742,0.608347270);
\draw[generalFenceInk!48,line width=.36pt] (1.147831871,-0.610313032) -- (0.362733795,0.931892802);
\draw[generalFenceInk!85,line width=.72pt,dash pattern=on 2.5pt off 1.5pt] (-1.163252762,0.192741612) -- (-1.102126983,-0.342322513) -- (-0.272659873,-0.962110489) arc[start angle=254.177393686,end angle=292.284862768,radius=1] -- (1.161495826,-0.249236506) -- (1.181073558,0.268114378) -- (0.362733795,0.931892802) arc[start angle=68.731817040,end angle=110.332074506,radius=1] -- (-1.163252762,0.192741612) -- cycle;
\draw[generalFenceInk,line width=.85pt] (-0.785776471,-0.099072287) -- (-0.272659873,-0.962110489) arc[start angle=254.177393686,end angle=292.284862768,radius=1] -- (0.853595939,-0.032331268) -- (0.362733795,0.931892802) arc[start angle=68.731817040,end angle=110.332074506,radius=1] -- (-0.785776471,-0.099072287) -- cycle;
\draw[black!80,line width=.65pt] (-1.211768873,0.617427080) -- (-1.457999121,-0.076410496) -- (-1.056706137,-0.739913603) -- (-0.272659873,-0.962110489) arc[start angle=254.177393686,end angle=292.284862768,radius=1] -- (1.147831871,-0.610313032) -- (1.479774589,0.025829562) -- (1.193948742,0.608347270) -- (0.362733795,0.931892802) arc[start angle=68.731817040,end angle=110.332074506,radius=1] -- (-1.211768873,0.617427080) -- cycle;
\draw[generalFenceArc,line width=1.15pt] (0.362733795,0.931892802) arc[start angle=68.731817040,end angle=110.332074506,radius=1];
\draw[generalFenceArc,line width=1.15pt] (-0.272659873,-0.962110489) arc[start angle=254.177393686,end angle=292.284862768,radius=1];
\end{scope}
\fill[black!85] (0.812048861,-0.604209901) circle[radius=.67pt];
\node[font=\scriptsize] at (0.908322839,-0.675843141) {$v_1$};
\fill[black!85] (1.111897307,0.025571266) circle[radius=.67pt];
\node[font=\scriptsize] at (1.231865585,0.028330280) {$v_2$};
\fill[black!85] (0.944129338,0.602263797) circle[radius=.67pt];
\node[font=\scriptsize] at (1.045298075,0.666799730) {$v_3$};
\fill[black!85] (-0.859432488,0.611252809) circle[radius=.67pt];
\node[font=\scriptsize] at (-0.957221709,0.680803282) {$v_4$};
\fill[black!85] (-1.099612180,-0.075646391) circle[radius=.67pt];
\node[font=\scriptsize] at (-1.219329230,-0.083882170) {$v_5$};
\fill[black!85] (-0.851722691,-0.732514467) circle[radius=.67pt];
\node[font=\scriptsize] at (-0.942703174,-0.810761203) {$v_6$};
\draw[black!75,fill=white,line width=.55pt] (0.346601770,0.922573874) circle[radius=.55pt];
\draw[black!75,fill=white,line width=.55pt] (-0.185579909,0.928317625) circle[radius=.55pt];
\draw[black!75,fill=white,line width=.55pt] (-0.281463744,-0.952489384) circle[radius=.55pt];
\draw[black!75,fill=white,line width=.55pt] (0.210383988,-0.916056837) circle[radius=.55pt];
\node at (0,-1.28) {(a) Convex projective domain};
\end{scope}
\begin{scope}[xshift=3.35cm]
\draw[black!18,line width=.5pt] (0,0) circle[radius=1];
\fill[black!3] (-0.766044443,0.642787610) arc[start angle=50.000000000,delta angle=-139.000000000,radius=0.373884679] arc[start angle=91.000000000,delta angle=-142.999998196,radius=0.334595320] arc[start angle=127.999997895,delta angle=-147.999997895,radius=0.286745386] arc[start angle=250.000000000,end angle=290.000000000,radius=1] arc[start angle=-160.000000000,delta angle=-150.000000000,radius=0.267949192] arc[start angle=-130.000000000,delta angle=-144.000000000,radius=0.324919696] arc[start angle=-94.000000000,delta angle=-143.999998142,radius=0.324919696] arc[start angle=-58.000002038,delta angle=-146.999997962,radius=0.296213495] arc[start angle=65.000000000,end angle=115.000000000,radius=1] arc[start angle=25.000000000,delta angle=-155.000000000,radius=0.221694663] -- cycle;
\fill[generalFenceFill] (-0.503132342,0.227361092) arc[start angle=14.657205562,delta angle=-31.329034354,radius=0.809784033] arc[start angle=46.045456978,delta angle=-66.045456978,radius=0.687280959] arc[start angle=250.000000000,end angle=290.000000000,radius=1] arc[start angle=-160.000000000,delta angle=-68.417506355,radius=0.649407593] arc[start angle=-168.439503924,delta angle=-30.154366805,radius=0.726542528] arc[start angle=-136.807545481,delta angle=-68.192454519,radius=0.687280959] arc[start angle=65.000000000,end angle=115.000000000,radius=1] arc[start angle=25.000000000,delta angle=-63.525872359,radius=0.649407593] -- cycle;
\fill[generalFenceDeep] (-0.347224053,0.051285488) arc[start angle=20.601747537,delta angle=-40.601747537,radius=1.428148007] arc[start angle=250.000000000,end angle=290.000000000,radius=1] arc[start angle=-160.000000000,delta angle=-42.634228295,radius=1.234897157] arc[start angle=-162.235191331,delta angle=-42.764808669,radius=1.303225373] arc[start angle=65.000000000,end angle=115.000000000,radius=1] arc[start angle=25.000000000,delta angle=-39.921605620,radius=1.257172299] -- cycle;
\draw[generalFenceInk!55,line width=.33pt,dash pattern=on 1.8pt off 1.4pt] (-0.422618262,0.906307787) arc[start angle=25.000000000,delta angle=-114.000000000,radius=0.649407593];
\draw[generalFenceInk!55,line width=.33pt,dash pattern=on 1.8pt off 1.4pt] (-0.766044443,0.642787610) arc[start angle=50.000000000,delta angle=-101.999999254,radius=0.809784033];
\draw[generalFenceInk!55,line width=.33pt,dash pattern=on 1.8pt off 1.4pt] (-0.999847695,-0.017452406) arc[start angle=91.000000000,delta angle=-111.000000000,radius=0.687280959];
\draw[generalFenceInk!55,line width=.33pt,dash pattern=on 1.8pt off 1.4pt] (0.342020143,-0.939692621) arc[start angle=-160.000000000,delta angle=-114.000000000,radius=0.649407593];
\draw[generalFenceInk!55,line width=.33pt,dash pattern=on 1.8pt off 1.4pt] (0.766044443,-0.642787610) arc[start angle=-130.000000000,delta angle=-107.999999169,radius=0.726542528];
\draw[generalFenceInk!55,line width=.33pt,dash pattern=on 1.8pt off 1.4pt] (0.997564050,-0.069756474) arc[start angle=-94.000000000,delta angle=-111.000000000,radius=0.687280959];
\draw[generalFenceInk!48,line width=.36pt] (-0.422618262,0.906307787) arc[start angle=25.000000000,delta angle=-76.999999520,radius=1.257172299];
\draw[generalFenceInk!48,line width=.36pt] (-0.766044443,0.642787610) arc[start angle=50.000000000,delta angle=-70.000000000,radius=1.428148007];
\draw[generalFenceInk!48,line width=.36pt] (0.342020143,-0.939692621) arc[start angle=-160.000000000,delta angle=-77.999999511,radius=1.234897157];
\draw[generalFenceInk!48,line width=.36pt] (0.766044443,-0.642787610) arc[start angle=-130.000000000,delta angle=-75.000000000,radius=1.303225373];
\draw[generalFenceInk!85,line width=.72pt,dash pattern=on 2.5pt off 1.5pt] (-0.503132342,0.227361092) arc[start angle=14.657205562,delta angle=-31.329034354,radius=0.809784033] arc[start angle=46.045456978,delta angle=-66.045456978,radius=0.687280959] arc[start angle=250.000000000,end angle=290.000000000,radius=1] arc[start angle=-160.000000000,delta angle=-68.417506355,radius=0.649407593] arc[start angle=-168.439503924,delta angle=-30.154366805,radius=0.726542528] arc[start angle=-136.807545481,delta angle=-68.192454519,radius=0.687280959] arc[start angle=65.000000000,end angle=115.000000000,radius=1] arc[start angle=25.000000000,delta angle=-63.525872359,radius=0.649407593] -- cycle;
\draw[generalFenceInk,line width=.85pt] (-0.347224053,0.051285488) arc[start angle=20.601747537,delta angle=-40.601747537,radius=1.428148007] arc[start angle=250.000000000,end angle=290.000000000,radius=1] arc[start angle=-160.000000000,delta angle=-42.634228295,radius=1.234897157] arc[start angle=-162.235191331,delta angle=-42.764808669,radius=1.303225373] arc[start angle=65.000000000,end angle=115.000000000,radius=1] arc[start angle=25.000000000,delta angle=-39.921605620,radius=1.257172299] -- cycle;
\draw[black!80,line width=.65pt] (-0.766044443,0.642787610) arc[start angle=50.000000000,delta angle=-139.000000000,radius=0.373884679] arc[start angle=91.000000000,delta angle=-142.999998196,radius=0.334595320] arc[start angle=127.999997895,delta angle=-147.999997895,radius=0.286745386] arc[start angle=250.000000000,end angle=290.000000000,radius=1] arc[start angle=-160.000000000,delta angle=-150.000000000,radius=0.267949192] arc[start angle=-130.000000000,delta angle=-144.000000000,radius=0.324919696] arc[start angle=-94.000000000,delta angle=-143.999998142,radius=0.324919696] arc[start angle=-58.000002038,delta angle=-146.999997962,radius=0.296213495] arc[start angle=65.000000000,end angle=115.000000000,radius=1] arc[start angle=25.000000000,delta angle=-155.000000000,radius=0.221694663] -- cycle;
\draw[generalFenceArc,line width=1.15pt] (0.422618262,0.906307787) arc[start angle=65.000000000,end angle=115.000000000,radius=1];
\draw[generalFenceArc,line width=1.15pt] (-0.342020143,-0.939692621) arc[start angle=250.000000000,end angle=290.000000000,radius=1];
\fill[black!85] (0.766044443,-0.642787610) circle[radius=.67pt];
\node[font=\scriptsize] at (0.857969776,-0.719922123) {$v_1$};
\fill[black!85] (0.997564050,-0.069756474) circle[radius=.67pt];
\node[font=\scriptsize] at (1.117271736,-0.078127251) {$v_2$};
\fill[black!85] (0.848048096,0.529919264) circle[radius=.67pt];
\node[font=\scriptsize] at (0.949813868,0.593509576) {$v_3$};
\fill[black!85] (-0.766044443,0.642787610) circle[radius=.67pt];
\node[font=\scriptsize] at (-0.857969776,0.719922123) {$v_4$};
\fill[black!85] (-0.999847695,-0.017452406) circle[radius=.67pt];
\node[font=\scriptsize] at (-1.119829419,-0.019546695) {$v_5$};
\fill[black!85] (-0.788010754,-0.615661475) circle[radius=.67pt];
\node[font=\scriptsize] at (-0.882572044,-0.689540852) {$v_6$};
\draw[black!75,fill=white,line width=.55pt] (0.422618262,0.906307787) circle[radius=.55pt];
\draw[black!75,fill=white,line width=.55pt] (-0.422618262,0.906307787) circle[radius=.55pt];
\draw[black!75,fill=white,line width=.55pt] (-0.342020143,-0.939692621) circle[radius=.55pt];
\draw[black!75,fill=white,line width=.55pt] (0.342020143,-0.939692621) circle[radius=.55pt];
\node at (0,-1.28) {(b) Ideal convex domain};
\end{scope}

\fill[black!3] (-.95,-1.61) rectangle (-.70,-1.51);
\draw[black!80,line width=.65pt] (-.95,-1.61) rectangle (-.70,-1.51);
\node[anchor=west] at (-.63,-1.56) {$\Omega$};
\fill[generalFenceFill] (.15,-1.61) rectangle (.40,-1.51);
\draw[generalFenceInk!85,line width=.72pt,dash pattern=on 2.5pt off 1.5pt]
 (.15,-1.61) rectangle (.40,-1.51);
\node[anchor=west] at (.47,-1.56) {$F_2$};
\fill[generalFenceDeep] (1.30,-1.61) rectangle (1.55,-1.51);
\draw[generalFenceInk,line width=.85pt] (1.30,-1.61) rectangle (1.55,-1.51);
\node[anchor=west] at (1.62,-1.56) {$F_3$};
\draw[generalFenceArc,line width=1.15pt] (2.45,-1.56) -- (2.70,-1.56);
\node[anchor=west,font=\footnotesize] at (2.77,-1.56) {continuous arcs};
\end{tikzpicture}
\caption{\small General-domain fences in a projective affine chart (a) and
in the Poincar\'e disk (b). In both figures, six isolated vertices form
two polygonal boundary chains separated by continuous arcs of extreme
points (orange); open circles mark the nonisolated arc endpoints.
The palest shading indicates $\Omega$, while the light and dark
blue shades depict the closures of $F_2$ and $F_3$, respectively,
with all defining cuts drawn. The continuous arcs
remain in both closures and produce no additional cuts.}
\label{fig:general-fence-arcs}
\end{figure}

An extreme point is \emph{isolated} if some neighborhood of it
in $\partial K$ contains no other extreme point.
The cuts remove blocks of at most $s-1$ isolated extreme points;
their endpoints need not be isolated. For $s=2$, they cut off
corners with a straight boundary segment on each side. If the
extreme-point set has no isolated points, then $F_s(\Omega)=\Omega$.
Arcs of nonisolated extreme points cannot be crossed by polygonal boundary chains,
so vertices separated by such an arc are not consecutive
(Figure~\ref{fig:general-fence-arcs}).

Infinite vertex sets occur naturally: for every integer $m\ge0$, the
differential $Q=z^m e^z\,dz^2$ gives a harmonic diffeomorphism
from $\C$ with complete $\partial$-energy metric \cite{WanAu}.
Its image has infinitely many isolated ideal vertices and exactly
one accumulation point \cite[Corollary~3.1]{ATW}.

Unlike the polygonal case, this fence need not be compactly contained
in $\Omega$. We therefore use its relative boundary
$\partial_\Omega F_s(\Omega)$ and set the distance to the empty set
equal to $+\infty$.

\begin{corollary}\label{cor:general-fences}
In the setting of Theorem~\ref{thm:pick}, set $k=m+3$ and
$d=d_\Omega$. In the setting of Theorem~\ref{thm:hopf}, set
$p=f(z_0)$, $k=m+2$, and $d=d_{\HH}$.
For every integer $s\ge2$ with $k\ge2s$, we have
$p\in F_s(\Omega)$ and
\[
 d\bigl(p,\partial_\Omega F_s(\Omega)\bigr)
 \ge\log\cot\frac{s\pi}{2k},
\]
with exactly the equality cases listed in
Corollary~\ref{cor:higher-fences}. In particular,
$p\in\Int F_s(\Omega)$ unless $k=2s$ and $\Omega$ is the
corresponding regular $2s$-gon; in that case $F_s(\Omega)=\{p\}$.
\end{corollary}

\begin{proof}
Use the same coordinates, centered radii $a,b$, and angle
$\vartheta=\arccos(a/b)$ as in the proof of
Corollary~\ref{cor:higher-fences}. In the Hopf case,
\eqref{eq:support-domain} identifies the image with the interior
of the convex hull of its finite ideal support directions, so its
centered circumradius in the Klein disk is $b=1$.
If $a=b$, the domain is a disk,
there are no cuts, and the conclusion is immediate. Otherwise,
Lemma~\ref{lem:chain-disk} applies to every defining boundary chain
of $j\le s$ segments. Its retained half-plane contains the closed
disk of radius $b\cos(j\vartheta)\ge b\cos(s\vartheta)$, so
\[
 \{x\in\R^2:|x|\le b\cos(s\vartheta)\}\subset F_s(\Omega).
\]
This common disk remains contained after intersecting any number
of cuts. The radial metric comparison in
\eqref{eq:fence-annulus-distance} gives the stated bound, since
$\vartheta\le\pi/k$, with strict inequality outside the corresponding
regular $k$-gon by Theorems~\ref{thm:pick} and~\ref{thm:hopf}.
The regular cases are already covered by
Corollary~\ref{cor:higher-fences}.
\end{proof}

\FloatBarrier
\section{Examples}\label{sec:examples}
The multiplicity conditions are uniform sufficient conditions, not
necessary ones. The following examples place all zeros, each simple,
inside a prescribed higher fence in both settings. We then construct
projective polygons with simple Pick zeros outside higher fences.

\begin{proposition}[Simple zeros inside higher fences]
\label{prop:simple-inside}
For every pair of integers $s\ge2$ and $n>2s$, the following hold.
\begin{enumerate}[label=\textup{(\roman*)},leftmargin=2em]
\item There is a convex projective $n$-gon $P$ whose $n-3$ Pick
zeros are distinct and simple and all belong to $\Int F_s(P)$.
\item There is an orientation-preserving harmonic diffeomorphism
$f:\C\to P$ onto an ideal $n$-gon, with complete
$\partial$-energy metric, whose $n-2$ Hopf zeros are distinct and
simple and whose images all belong to $\Int F_s(P)$.
\end{enumerate}
\end{proposition}

\begin{proof}
For \textup{(i)}, put $d=n-3$ and consider
\[
 C_\varepsilon=(z^d-\varepsilon)\,dz^3,
 \qquad \varepsilon\ne0,\qquad \varepsilon\longrightarrow0.
\]
The limiting differential $C_0=z^d\,dz^3$ corresponds to a regular
$n$-gon $P_0$, with the image $p_0$ of $0$ at its center.
Since $n>2s$, every defining diagonal of $F_s(P_0)$ has angular
span $2\pi s/n<\pi$, so $p_0\in\Int F_s(P_0)$.
Choose the affine-sphere reconstructions with continuously varying
initial frames. The continuity of the complete metrics and of
the reconstructed affine spheres and their polygons
\cite[Theorem~5.4 and the proof of Theorem~8.1]{DW} gives local
uniform convergence of their projectivizations
$f_\varepsilon:\C\to P_\varepsilon$ to $f_0$, together with
convergence of the vertices of $P_\varepsilon$ to those of $P_0$.
Thus a fixed neighborhood of $p_0$ lies in $\Int F_s(P_\varepsilon)$
for small $\varepsilon$.
All $d$ roots of $z^d-\varepsilon$ are simple and converge to $0$,
so their images lie in that neighborhood. These are all the Pick
zeros of $P_\varepsilon$.

For \textup{(ii)}, put $d=n-2$ and take
\[
 Q_\varepsilon=q_\varepsilon(z)\,dz^2
 =(z^d-\varepsilon)\,dz^2.
\]
Let $f_\varepsilon:\C\to P_\varepsilon$ be the corresponding
harmonic diffeomorphisms with complete $\partial$-energy metrics
$I_{f_\varepsilon}=H_\varepsilon|dz|^2$
\cite[Theorem~1.1 and \S2]{HTTW}.
Their Gauss equation is
\[
 \Delta\log H_\varepsilon
 =2\left(H_\varepsilon-
          \frac{|q_\varepsilon|^2}{H_\varepsilon}\right).
\]
This is the order-two case of the vortex equation in
\cite[Theorem~5.4]{DW}. The metrics are complete and negatively
curved, since $|q_\varepsilon|<H_\varepsilon$.
That theorem gives local $C^1$ continuity
of $\log H_\varepsilon$ in the polynomial coefficients.
The CMC structure equations, with fixed initial oriented orthonormal
frames, therefore give local uniform convergence
$f_\varepsilon\to f_0$. Normalize the target so that
$f_\varepsilon(0)=0$ in the Klein disk. Uniqueness and cyclic
symmetry identify $P_0$ with the regular ideal $n$-gon centered at $0$.

Choose $a$ such that
\[
 \cos\frac{\pi}{2s}<a<\cos\frac{\pi}{n}.
\]
The closed Euclidean disk $\overline B_a(0)$ is compactly contained
in $P_0$. Let $b_1,\ldots,b_n$ be the ideal vertices of $P_0$, and choose
$t<1$ so close to $1$ that $t\cos(\pi/n)>a$.
The points $tb_i$ lie in $P_0$, and their convex hull contains
$\overline B_a(0)$ in its interior. Local convergence at the fixed
points $f_0^{-1}(tb_i)$, together with convexity of $P_\varepsilon$, gives
$\overline B_a(0)\subset P_\varepsilon$ for small $\varepsilon$.
Applying Lemma~\ref{lem:chain-disk} to each chain of $s$ sides,
with $b=1$, yields
\[
 \overline B_r(0)\subset F_s(P_\varepsilon),
 \qquad r=\cos\bigl(s\arccos a\bigr)>0.
\]
All $d$ roots of $q_\varepsilon$ are simple when $\varepsilon\ne0$
and converge to $0$. Local uniform convergence of the maps implies
that their images converge to $0$, so they all lie in
$B_r(0)\subset\Int F_s(P_\varepsilon)$ for small $\varepsilon$.
\end{proof}

\begin{remark}
Let $k\ge3$, $s\ge2$, and let $P$ be a polynomial of degree $s-1$
with $P(0)\ne0$. Consider
\[
 Q=z^{k-2}P(z^k)\,dz^2.
\]
The corresponding complete harmonic diffeomorphism
$f:\C\to\Omega$ has an ideal $ks$-gon as its image.
For $s=2$ and $P(w)=w-c$, this is the family
$(z^{2m}-cz^{m-1})\,dz^2$ with $m=k-1$ studied in \cite{AuWanSymmetry}.
Only cyclic symmetry is needed for the following fence calculation.

Put $p=f(0)$. Invariance of $Q$ under
$z\mapsto e^{2\pi i/k}z$, uniqueness of the complete metric
\cite[Theorem~3.1]{LiVortex}, and geometric reconstruction give a
target rotation fixing $p$. Since $Q(0)=0$, the derivative of $f$
there is conformal, so the rotation has angle $2\pi/k$.
In the Klein disk centered at $p$, it shifts the cyclic vertex
labels by $s$ positions. Every defining chord of $F_s(\Omega)$
therefore has angular span $2\pi/k$ and Euclidean support distance
$\cos(\pi/k)$ from $p$; the original sides lie farther away.
Consequently,
\[
 d_{\HH}\bigl(p,\partial F_s(\Omega)\bigr)
 =\operatorname{arctanh}\!\left(\cos\frac{\pi}{k}\right)
 =\log\cot\frac{\pi}{2k}.
\]
If the roots of $P$ are distinct and nonzero, the zero at $0$ has
order $k-2$ and all other zeros are simple. In particular, $k=3$
gives $3s-2$ simple zeros and the exact distance $\frac12\log3$.

To place all zero images inside the fence, take
$P_a(w)=w^{s-1}+a$ with $a\ne0$ tending to $0$, and write
$f_a:\C\to\Omega_a$ and $p_a=f_a(0)$.
The metric and reconstruction continuity used in the preceding
proof applies to this family. All zeros converge to $0$, so their
images approach $p_a$ in hyperbolic distance. The displayed formula
provides a fixed disk centered at $p_a$ inside $F_s(\Omega_a)$;
hence all zero images lie in its interior for sufficiently small
$a\ne0$. When $k=3$, this recovers the family
$(z^{3s-2}+az)\,dz^2$ with all zeros simple.
\end{remark}

\begin{proposition}[Simple Pick zeros outside higher fences]
\label{prop:simple-pick-outside}
For every integer $s\ge3$ and every integer $n\ge2s$, there is a
convex projective $n$-gon with at least $2s-4$ distinct simple Pick
zeros outside $F_s(P)$.
\end{proposition}
\begin{proof}
Start with the regular $(2s-1)$-gon $P_0$. Its central Pick zero
has order $2s-4$. The diagonal joining $v_0$ to $v_s$ has forward
angular span $2\pi s/(2s-1)>\pi$. Its retained half-plane therefore
excludes a small closed disk $K$ about the center, containing no
other zero and no zero on its boundary.

Replace an edge on the complementary arc by a shallow outward
strictly convex chain with $n-(2s-1)$ additional vertices.
This gives $n$-gons $P_j\to P_0$ whose selected diagonal still
excludes $K$. For fixed $n$, simple-root polynomials are dense in
the degree-$(n-3)$ polynomial moduli space. The Dumas--Wolf
homeomorphism \cite[Theorem~A]{DW} therefore allows arbitrarily
small perturbations of the $P_j$ for which all Pick zeros are
simple. Fixing four original vertices as a projective normalization,
choose these perturbations to preserve convergence and the strict
exclusion of $K$.

The canonical potentials, Blaschke metrics and Pick tensors
converge on compact subsets by
\cite[Theorem~4.4 and Corollary~4.5]{DW}.
Choose convergent oriented orthonormal coframes extending smoothly
over the entire disk $K$ and write the complex Pick tensors as
$c_j(e_j^1+ie_j^2)^3$, with $j=0$ denoting the limiting data.
Uniform convergence and $c_0\ne0$ on $\partial K$ give
$|c_j-c_0|<|c_0|$ there for all large $j$.
Thus the straight-line homotopy between $c_0$ and $c_j$ avoids $0$
on $\partial K$, so their boundary winding numbers agree.
The change from each such coframe to a local holomorphic frame is
a smooth nonvanishing factor, so it contributes no local zero index.
The argument principle therefore identifies each local index with
the holomorphic vanishing order, and the sum of these indices with
the boundary winding number. Hence the total Pick zero order in
$K$ remains $2s-4$. All these zeros are simple and $K$ is disjoint
from $F_s(P_j)$, proving the claim.
\end{proof}

\end{document}